\documentclass[12pt]{amsart}
\usepackage{amsmath}
\usepackage{amsfonts}
\usepackage{amssymb}
\usepackage[all]{xy}
\usepackage{bbding}
\usepackage{txfonts}
\usepackage{amscd}
\usepackage{xspace}
\usepackage[T1]{fontenc}
\usepackage[shortlabels]{enumitem}
\usepackage{ifpdf}
\ifpdf
  \usepackage[colorlinks,final,backref=page,hyperindex]{hyperref}
\else
  \usepackage[colorlinks,final,backref=page,hyperindex]{hyperref}
\fi
\usepackage{tikz}
\usepackage[active]{srcltx}

\newtheorem{thm}{Theorem}[section]
\newtheorem{lem}[thm]{Lemma}
\newtheorem{cor}[thm]{Corollary}
\newtheorem{pro}[thm]{Proposition}
\newtheorem{ex}[thm]{Example}
\newtheorem{rmk}[thm]{Remark}
\newtheorem{defi}[thm]{Definition}

\newcommand {\emptycomment}[1]{}
\newcommand{\lon }{\,\rightarrow\,}
\newcommand{\be }{\begin{equation}}
\newcommand{\ee }{\end{equation}}

\newcommand{\pf}{\noindent{\bf Proof.}\ }

\newcommand{\Real}{\mathbb R}
\newcommand{\Comp}{\mathbb C}
\newcommand{\Nat}{\mathbb N}
\newcommand{\Integ}{\mathbb Z}

\newcommand{\huaL}{\mathcal{L}}
\newcommand{\huaR}{\mathcal{R}}

\newcommand{\huaF}{\mathcal{F}}

\newcommand{\huaM}{\mathcal{M}}

\newcommand{\huaD}{\mathcal{D}}

\newcommand{\huaO}{\mathcal{O}}

\def\qed{\hfill ~\vrule height6pt width6pt depth0pt}

\newcommand{\half}{\frac{1}{2}}

\newcommand{\Id}{\rm{Id}}

\newcommand{\br}[1]{   [ \cdot,    \cdot  ]   }

\newcommand{\dM}{\mathrm{d}}

\newcommand{\ad}{\mathrm{ad}}

\newcommand{\CE}{\mathrm{CE}}

\newcommand{\FGV}{\mathrm{FGV}}

\newcommand{\GD}{\mathrm{GD}}
\newcommand{\PGD}{\mathrm{PGD}}
\allowdisplaybreaks[3]

\begin{document}

\title[Poisson conformal  algebras and deformation quantization]{Cohomology and deformation quantization of  Poisson conformal  algebras}

\author{Jiefeng Liu}
\address{School of Mathematics and Statistics, Northeast Normal University, Changchun 130024, China}
\email{liujf534@nenu.edu.cn}

\author{Hongyu Zhou}
\address{School of Mathematics and Statistics, Northeast Normal University, Changchun 130024, China}
\email{zhouhy929@nenu.edu.cn}

\begin{abstract}
In this paper, we first recall the notion of (noncommutative)  Poisson conformal  algebras and describe some constructions of them. Then we study the formal distribution (noncommutative) Poisson algebras and  coefficient (noncommutative) Poisson algebras.  Next, we introduce the notion of conformal formal deformations of commutative associative conformal algebras and show that Poisson conformal  algebras are the corresponding semi-classical limits. At last, we develop the cohomology theory of noncommutative Poisson conformal  algebras and use this cohomology to study their deformations.
\end{abstract}



\maketitle

\tableofcontents

\section{Introduction}
Poisson algebras were originally coming from the study of the Hamiltonian mechanics and then appeared in many fields in mathematics and mathematical physics, such as Poisson geometry (\cite{Wei77,Vaisman1}), classical and quantum mechanics (\cite{Arn78, BFFLS}),  algebraic geometry (\cite{GK04,Pol97}), deformation theory and  quantization (\cite{BFFLS1,BFFLS2,Kon03}) and quantum groups (\cite{CP1,Dr87}).

The notion of a Lie conformal algebra (also called a vertex Lie algebra  in \cite{Prim}) was introduced by Kac as an axiomatic description of the operator product expansion  of chiral fields in two-dimensional conformal field theory (\cite{Kac}). It turns out to be a valuable tool in the study of vertex algebras (\cite{BDK2,Bakalov,Kac}), Poisson vertex algebras (also called vertex Poisson algebras) (\cite{BDK,LHS}) and Hamiltonian formalism in the theory of nonlinear evolution equations (\cite{BSK}). The structure theory (\cite{DAK}), representation theory (\cite{CK}) and cohomology theory (\cite{BKV,DK09}) of Lie conformal algebras have been well developed. Associative conformal algebras naturally come from studying the  representation theory of Lie conformal algebras. Moreover, any  Lie conformal
 algebras appeared in physics can be embedded into an associative conformal algebra, which is different form the Lie algebra case in general (see \cite{Rot1,Rot2} for more details).  Therefore, the study of associative conformal
algebras provides some information on the structures and representations of Lie
conformal algebras. See \cite{BKV,Dolg,Kac97,Kol06,Kol19} for more details on associative conformal algebras. 

Let $H$ be  a Hopf algebra generated by a primitive element $\partial$. The class $\huaM^*(H)$ of $H$-modules introduced in \cite{BAK} forms a pseudo-tensor category. Then a Lie (or associative ) conformal algebra can be seen as a morphism of operads from the operad Lie (or Ass) to $\huaM^*(H)$. In particular, for a morphism  from the operad Pois to $\huaM^*(H)$, the author in \cite{Kol20} introduced the notion of a Poisson conformal algebra, which consists of a Lie conformal algebra and a commutative associative conformal algebra satisfying a  conformal Leibniz rule. As  in the case of Poisson algebras, Poisson conformal algebras are also constructed by associated graded conformal algebras ${{\rm gr} ~U}$ of associative conformal envelopes of Lie conformal algebras. Based on the close relations between Poisson conformal algebras and representations of Lie conformal algebras, Poisson conformal algebras play an important role in the study of Ado-type problems for Lie conformal algebras (\cite{KKP}). It is well-known that Poisson algebras can be understood as semi-classical limits of formal deformations of
commutative associative algebras. So it is natural to ask whether Poisson conformal algebras can be understood as certain semi-classical limits.

The first aim of this paper is to answer the above question. Because of the locality condition for associative conformal algebras, for a vector space $A$, the $\Comp[[\hbar]]$-module $A[[\hbar]]$ does not carry the structure of an associative conformal algebra in general. In order to solve this problem, we introduce the notion of an $\hbar$-adic associative conformal algebra, which can be thought of as an inverse limit of the Lie conformal algebras $A[[\hbar]] /\hbar^n A[[\hbar]]$ over $ \mathbb{\Comp}[[\hbar]]$. On this basis,  we introduce the notion of conformal formal deformation of commutative associative conformal algebras and show that Poisson conformal algebras are the corresponding semi-classical limits. This result is parallel to that the semi-classical limit of an associative formal deformation of a commutative associative algebra is a Poisson algebra and the semi-classical limit of a vertex formal deformation of a commutative vertex algebra is a  Poisson vertex algebra (\cite{FB,LHS}).  Furthermore,  we show that conformal infinitesimal deformation and extension of conformal $n$-deformation to conformal $(n+1)$-deformation of a commutative associative conformal algebra are classified by the second and the third cohomology groups of the associative conformal algebra respectively.

The notion of a noncommutative Poisson algebra was first given by Xu in \cite{Xu}, which is especially suitable for geometric situation. In \cite{FGV}, Flato, Gerstenhaber and Voronov introduced a more general notion of a Leibniz pair and study its cohomology and deformation theory. In particular, they gave the cohomology theory of a noncommutative Poisson algebra associated to a module using an innovative bicomplex.  Recently, Bao and Ye developed the cohomology theory of noncommutative Poisson algebras associated to  quasi-modules through Yoneda-Ext groups and projective resolutions in \cite{CohomologyPA1,CohomologyPA2}. Noncommutative Poisson algebras had been studied by many authors from different aspects \cite{Kubo1,Kubo3,Kubo4,YYY}. A Poisson algebra in the usual sense is the one where the associative multiplication on $P$ is commutative. In this paper, we develop  the cohomology theory of noncommutative Poisson conformal algebras. Here the noncommutative Poisson conformal algebra means that the associative $\lambda$-multiplication on Poisson conformal algebra $P$ is noncommutative. We find that this cohomology can control simultaneous deformations of the associative and Lie conformal $\lambda$-multiplications in a noncommutative Poisson conformal algebra. Note that there exists another cohomology theory of Poisson conformal algebras, which is similar to the one of Poisson algebras  given by Lichnerowicz in \cite{Lic77}. We will study this cohomology in the further.

The paper is organized as follows. In Section \ref{sec:construction}, we recall the definition of a (noncommutative) Poisson conformal algebra and give some constructions. In particular, we introduce the notion of a  quadratic (noncommutative) Poisson conformal algebra, which can be constructed by (noncommutative) {\rm PGD}-algebras.  In Section \ref{sec: distribution}, we study formal  distribution (noncommutative)  Poisson algebras and coefficient (noncommutative) Poisson algebras. The later are useful to study infinite-dimensional (noncommutative)  Poisson algebras. In Section \ref{sec:deformation quantization}, first we introduce the notion of an $\hbar$-adic associative conformal algebra.  Then we introduce the notion of conformal formal deformations of commutative associative conformal algebras
and show that Poisson conformal algebras are the corresponding
semi-classical limits. Furthermore, we study conformal infinitesimal
deformations and extensions of conformal $n$-deformations to conformal
$(n+1)$-deformations of a commutative associative conformal algebra. In Section \ref{sec:cohomology}, we develop a cohomology theory of noncommutative Poisson conformal algebras with coefficients in an arbitrary module. This cohomology is also a bicomplex, which is induced by the Lie conformal algebra cohomology and the associative conformal algebra cohomology. In particular, if we replace the  associative conformal algebra cohomology by the commutative associative conformal algebra cohomology in the bicomplex, we can obtain the cohomology  of  Poisson conformal algebras. Furthermore, we use cohomology of  noncommutative Poisson conformal algebras to study their linear deformations.

\section{(Noncommutative) Poisson conformal algebras and some constructions}\label{sec:construction}

\begin{defi}
	A {\bf conformal algebra} $A$ is a $ \mathbb{C}[\partial] $-module endowed with a $ \mathbb{C} $-bilinear map $ A\times A\rightarrow A[\lambda] $, denoted by $ a\times b\mapsto a_{\lambda}b $, satisfying
	\begin{eqnarray}\label{eq:conformal seq}
		(\partial a)_{\lambda}b=-\lambda a_{\lambda}b,\quad a_{\lambda}(\partial b)=(\partial+\lambda)a_{\lambda}b.
	\end{eqnarray}
\end{defi}

\begin{defi}
An {\bf associative conformal algebra} $A$ is a $ \mathbb{C}[\partial] $-module endowed with a $ \mathbb{C} $-bilinear map $ A\times A\rightarrow A[\lambda] $, denoted by $ a\times b\mapsto a\circ_{\lambda}b $ such that $A$ is a conformal algebra and  satisfies
\begin{eqnarray}
\label{eq:aca1}a\circ _{\lambda}(b\circ_{\mu}c)&=&(a\circ_{\lambda}b)\circ_{\lambda+\mu}c,\quad\forall~ a,b,c\in A.
\end{eqnarray}
An associative conformal algebra is called {\bf finite} if it is finitely generated as a $\Comp[\partial]$-module.
\end{defi}

An associative conformal algebra $(A,\circ_\lambda)$ is called {\bf commutative} if it satisfies
$${a\circ_{\lambda}b}=b\circ_{-\lambda-\partial}a,\quad \forall~a,b\in A.$$

\emptycomment{Define the $n$th product on an associative conformal algebra $A$ for $n\in \Nat$ by
\begin{eqnarray*}
  a\circ_{\lambda}b=\sum_{n=0}^{\infty}\lambda^{(n)}(a_{(n)}b).
\end{eqnarray*}
 All the conditions of associative conformal algebra can be translated in terms of $n$th products as follows:
\begin{eqnarray}
  (\partial a)_{(n)}b&=&-na_{(n-1)}b,\quad (a)_{(n)}(\partial b)=\partial(a_{(n-1)}b)+ja_{(n-1)}b,\\
a_{(m)}(b_{(n)} c)&=&\sum_{j=0}^{m}\binom{m}{j}(a_{(j)}b)_{(m+n-j)}c.
\end{eqnarray}
Furthermore, the commutative condition in a commutative associative algebra can be written by
\begin{equation}
a_{(n)}b=\sum_{j=0}^{\infty}(-1)^{j+n}\partial^{(j)}(b_{(n+j)}a).
\end{equation}}

\begin{defi}
A {\bf Lie conformal algebra} $A$ is a $ \mathbb{C}[\partial] $-module endowed with a $ \mathbb{C} $-bilinear map $ A\times A\rightarrow A[\lambda] $, denoted by $ a\times b\mapsto [a_{\lambda}b] $ such that $A$ is a conformal algebra and  satisfies
\begin{eqnarray}
\label{eq:Lie conformal2}  {[a_{\lambda}b]}&=&-[b_{-\lambda-\partial}a],\\
 \label{eq:Lie conformal3} {[a_{\lambda}[b_{\mu}c]]}&=&[[a_{\lambda}b]_{\lambda+\mu}c]+[b_{\mu}[a_{\lambda}c]],\quad\forall~ a,b,c\in A.
\end{eqnarray}
A Lie conformal algebra is called {\bf finite} if it is finitely generated as a $\Comp[\partial]$-module.
\end{defi}

\emptycomment{Define the $n$th product on a Lie conformal algebra $A$ for $n\in \Nat$ by
\begin{eqnarray*}
	[a_{\lambda}b]=\sum_{n=0}^{\infty}\lambda^{(n)}(a_{[n]}b),
\end{eqnarray*}
where $\lambda^{(n)}=\lambda^n/n!$ and $a,b\in A$. Note that the formal power series $[a_{\lambda}b]$ is well-defined, then for large enough $n$, $a_{[n]}b=0$. All the conditions of Lie conformal algebra can be translated in terms of $n$th products as follows:
\begin{eqnarray}
	(\partial a)_{[n]}b&=&-na_{[n-1]}b,\quad
	a_{[n]}b=-\sum_{j=0}^{\infty}(-1)^{j+n}\partial^{(j)}(b_{[n+j]}a),\\
	a_{[m]}(b_{[n]}c)&=&\sum_{j=0}^{m}\binom{m}{j}(a_{[j]}b)_{[m+n-j]}c+b_{[n]}(a_{[m]}c),
\end{eqnarray}
where $ a,b,c\in A$.}

Let $(A,\circ_{\lambda})$ be an associative conformal algebra. It is not hard to see that the $\lambda$-bracket
	\begin{equation}\label{eq:commmutative bracket}
		[a_\lambda b]_c={a\circ_{\lambda}b}-b\circ_{-\lambda-\partial}a,\quad \forall~a,b\in A
	\end{equation}
	defines a Lie conformal algebra on $A$.

\begin{defi}
 A {\bf  Poisson conformal algebra} $P$ is a $ \mathbb{C}[\partial] $-module endowed with a $ \mathbb{C} $-bilinear map $ P\times P\rightarrow P[\lambda] $, denoted by $ a\times b\mapsto a\circ_{\lambda}b $ and a $ \mathbb{C} $-bilinear map $ P\times P\rightarrow P[\lambda] $, denoted by $ a\times b\mapsto [a_{\lambda}b] $ such that  $(P,\circ_\lambda)$ is a commutative associative conformal algebra, $(P,[\cdot_\lambda\cdot])$ is a Lie conformal algebra and the following Leibniz rule holds:
\begin{eqnarray}\label{eq:Poisson conformal Leibniz}	[a_{\lambda}(b\circ_{\mu}c)]=[a_{\lambda}b]\circ_{\lambda+\mu}c+b\circ _{\mu}[a_{\lambda}c],\quad \forall~a,b,c\in P.
	\end{eqnarray}
If the associative $\lambda$-multiplication in the Poisson algebra $P$ is noncommutative, we call it a  {\bf noncommutative Poisson conformal algebra}.
\end{defi}
Let $(P,\circ_{\lambda},[\cdot_\lambda\cdot])$ be a (noncommutative) Poisson conformal algebra. Define the $n$-th products on the associative conformal algebra $(P,\circ_{\lambda})$ and the Lie conformal algebra $(P,[\cdot_\lambda\cdot])$   for $n\in \Nat$  respectively by
\begin{eqnarray*}
	a\circ_{\lambda}b=\sum_{n=0}^{\infty}\lambda^{(n)}(a_{(n)}b),\quad
	{	[a_{\lambda}b]}=\sum_{n=0}^{\infty}\lambda^{(n)}(a_{[n]}b),
\end{eqnarray*}
where $\lambda^{(n)}=\lambda^n/n!$ and $a,b\in P$. Then the conformal Leibniz rule  of (noncommutative) Poisson conformal algebras can be translated in terms of $n$-th products of Lie conformal algebras and associative conformal algebras as follows:
	\begin{eqnarray}\label{eq:conformal Leibniz}
	a_{[m]}(b_{(n)}c)=\sum_{j=0}^{m}\binom{m}{j}(a_{[j]}b)_{(m+n-j)}c+b_{(n)}(a_{[m]}c),\quad\forall~a,b,c\in P.
\end{eqnarray}

\begin{defi}
Let $(P_1,\circ_{1\lambda},[\cdot_\lambda\cdot]_1,\partial_1)$ and $(P_2,\circ_{2\lambda},[\cdot_\lambda\cdot]_2,\partial_2)$ be two (noncommutative) Poisson conformal algebras. A $\Comp$-linear map $\phi$ is called a {\bf homomorphism} from $P_1$ to $P_2$ if $\phi$ satisfies
\begin{eqnarray*}
\phi(\partial_1 a)&=&\partial_2\phi(a),\\
\phi(a\circ_{1\lambda}b)&=&\phi(a)\circ_{2\lambda}\phi(b),\\
\phi([a_{\lambda}b]_1)&=&[\phi(a)_{\lambda}\phi(b)]_2,\quad\forall~a,b\in P_1.
\end{eqnarray*}
\end{defi}

\begin{ex}\label{ex:direct sum of HMA}{
  Let $(P_1,\circ_{1\lambda},[\cdot_\lambda\cdot]_1,\partial_1)$ and $(P_2,\circ_{2\lambda},[\cdot_\lambda\cdot]_2,\partial_2)$ be two (noncommutative) Poisson conformal algebras. Then $(P_1\oplus P_2,\circ_{\oplus\lambda},[\cdot_\lambda\cdot]_{\oplus},\partial=\partial_1\oplus \partial_2)$ is a (noncommutative) Poisson conformal algebra, where the $\lambda$-product $\circ_{\oplus\lambda}$ and the $\lambda$-bracket $[\cdot_\lambda\cdot]_{\oplus}$ are given by
  \begin{eqnarray*}
   ( a_1+ a_2) \circ_{\oplus\lambda} (b_1+ b_2)&=& a_1\circ_{1\lambda} b_1+a_2\circ_{2\lambda} b_2,\\
   {[  ( a_1+ a_2)_{\lambda} (b_1+ b_2)]_{\oplus }}&=&[{a_1}_\lambda b_1]_1+ [{a_2}_\lambda b_2]_2,
  \end{eqnarray*}
  for all $a_1,b_1\in P_1,a_2,b_2\in P_2.$}
\end{ex}

\begin{ex}
	Let $( A,\circ,[\cdot,\cdot] )$ be an ordinary (noncommutative) Poisson algebra. Then $ P=\mathbb{C}[\partial]\otimes A $ equipped with $\lambda$-operations
$$ a\circ_{\lambda}b=a\circ b,\quad [a_{\lambda}b]=[a,b], \quad\forall~ a,b\in A $$
is a (noncommutative) Poisson conformal algebra.
\end{ex}

\begin{ex}
  Let $(A,\circ_\lambda)$ be an associative conformal algebra. Then $(A,\circ_\lambda,[\cdot_\lambda \cdot]_c)$ is a noncommutative Poisson conformal algebra, where the $\lambda$-bracket $[\cdot_\lambda \cdot]_c$ is given by \eqref{eq:commmutative bracket}.
\end{ex}

In the following, we will give the notion of  quadratic (noncommutative) Poisson conformal algebras.  First, we recall the notion of quadratic Lie conformal algebras.
\begin{defi}
	If a Lie conformal algebra $A=\Comp[\partial] V$ with a free $\Comp[\partial]$-module and the $\lambda$-bracket is of the following form:
	$$[a_\lambda b]=\partial u+\lambda v+w,\quad \forall~a,b\in V,$$
	where $u,v,w\in V$, then $A$ is  called a {\bf quadratic Lie conformal algebra}.
\end{defi}

Recall that a {\bf Novikov algebra} is a pair $(A,\ast)$, where $A$ is a vector space and $\ast:A\otimes A\rightarrow A$ is a bilinear multiplication satisfying that for all $a,b,c\in A$,
\begin{eqnarray*}
	\label{eq:Nov1}  (a\ast b)\ast c&=&(a\ast c)\ast b,\\
	\label{eq:Nov2}  (a\ast b)\ast c-a\ast (b\ast c)&=&(b\ast a)\ast c-b\ast(a\ast c).
\end{eqnarray*}

Recall that a {\bf Gel'fand-Dorfman bialgebra (or {\rm GD}-algebra)} is a triple $(A,\ast,[\cdot,\cdot])$ such that $(A,[\cdot,\cdot])$ is a Lie algebra, $(A,\ast)$ is a Novikov algebra and they satisfy the following condition:
\begin{equation*}\label{eq:GD condition}
	[a\ast b,c]+[a,b]\ast c-a\ast[b,c]-[a\ast c,b]-[a,c]\ast b=0,\quad\forall~a,b,c\in A.
\end{equation*}

\begin{pro}{\rm(\cite{Gel,XP})}\label{pro:quadratic Lie conformal algebra}
	A Lie conformal algebra $A=\Comp[\partial] V$ is quadratic if and only if $V$ is a {\rm GD}-algebra, where the correspondence is given as follows:
	$$[a_\lambda b]=\partial(b\ast a)+\lambda(a\bullet b)+[b,a],$$
	where $a\bullet b=a\ast b+b\ast a$ for $a,b\in V$.
\end{pro}

\begin{defi}
  If the associative conformal algebra in  a (noncommutative) Poisson conformal algebra $(P=\Comp[\partial]\otimes V,\circ_\lambda,[\cdot_\lambda\cdot])$ is a current associative conformal algebra and the Lie conformal algebra in $P$ is a quadratic Lie conformal algebra,  then $P$ is called a {\bf quadratic (noncommutative) Poisson conformal algebra}.
\end{defi}

\begin{pro}\label{pro:quadratic comformal Poisson}
  A (noncommutative) Poisson conformal algebra $P=\Comp[\partial] V$ is quadratic if and only if $(V,\ast,[\cdot,\cdot])$ is a {\rm GD}-algebra and $(V,\circ,[\cdot,\cdot])$ is a  (noncommutative) Poisson algebra such that the following equations hold:
  \begin{eqnarray}
\label{eq:quad Poi1}    (b\circ c)\ast a&=&b\circ(c\ast a)~=~(b\ast a)\circ c,\\
 \label{eq:quad Poi2}   a\ast(b\circ c)&=&(a\ast b)\circ c+b\circ (a\ast c),\quad\forall~a,b,c\in V.
  \end{eqnarray}
\end{pro}
\begin{proof}
  By Proposition \ref{pro:quadratic Lie conformal algebra}, $(P,[\cdot_\lambda\cdot])$ is a quadratic Lie conformal algebra if and only if $(V,\ast,[\cdot,\cdot])$ is a {\rm GD}-algebra. By the Leibniz rule of the (noncommutative) Poisson conformal algebra, we have
  \begin{eqnarray*}
    &&[a_{\lambda}(b\circ_{\mu}c)]-[a_{\lambda}b]\circ_{\lambda+\mu}c-b\circ _{\mu}[a_{\lambda}c]\\
    &=&\partial\big((b\circ c)\ast a-b\circ(c\ast a)\big)+\mu\big((b\ast a)\circ c-b\circ(c\ast a)\big)+\lambda\big(a\ast(b\circ c)\\
    &&+(b\circ c)\ast a+(b\ast a)\circ c-(a\ast b)\circ c-(b\ast a)\circ c-b\circ (a\ast c)-b\circ(c\ast a)\big)\\
    &&-([a,b\circ c]-[a,b]\circ c-b\circ [a,c]).
  \end{eqnarray*}
  Thus the  Leibniz rule holds if and only if
  \begin{eqnarray}
 \label{eq:quad proof Poi1}   (b\circ c)\ast a&=&b\circ(c\ast a),\quad (b\ast a)\circ c=b\circ(c\ast a),\\
   \nonumber a\ast(b\circ c)&+&(b\circ c)\ast a+(b\ast a)\circ c-(a\ast b)\circ c\\
\label{eq:quad proof Poi2}      &&-(b\ast a)\circ c)-b\circ (a\ast c)-b\circ(c\ast a)=0,\\
\label{eq:quad proof Poi3}    { [a,b\circ c]}&=&[a,b]\circ c+b\circ [a,c].
  \end{eqnarray}
 By \eqref{eq:quad proof Poi1},  \eqref{eq:quad Poi1} holds. Taking \eqref{eq:quad proof Poi1} into \eqref{eq:quad proof Poi2}, \eqref{eq:quad Poi2} follows. By \eqref{eq:quad proof Poi3}, $(V,\circ,[\cdot,\cdot])$ is a  (noncommutative) Poisson algebra.
 
 The converse can be proved similarly.
\end{proof}

\begin{defi}
  A {\bf (noncommutative) {\rm PGD}-algebra} is a quadruple $(A,\ast,\circ,[\cdot,\cdot])$ such that $(A,\ast,[\cdot,\cdot])$ is a {\rm GD}-algebra and $(A,\circ,[\cdot,\cdot])$ is a (noncommutative) Poisson algebra satisfying \eqref{eq:quad Poi1} and \eqref{eq:quad Poi2}.
\end{defi}

The following quadratic Poisson conformal algebra was given in \cite{Kol20}. By the relation between the {\rm PGD}-algebras and quadratic Poisson conformal algebras, we give a  easier proof here.
\begin{pro}\label{ex:derivation-HMA}
 Let $(A,\circ,[\cdot,\cdot])$ be a Poisson algebra with a derivation $D$. The operation
  \begin{eqnarray}\label{eq:Novikov product}
    a\ast b&=&a\circ D (b),\quad\forall~a,b\in A
  \end{eqnarray}
 makes $(A,\ast,\circ,[\cdot,\cdot])$ being a $\PGD$-algebra. Furthermore, $(P,\circ_\lambda,[\cdot_\lambda\cdot])$ is a Poisson conformal algebra, where $\circ_\lambda$ and $[\cdot_\lambda\cdot]$ are given by
  \begin{eqnarray}
  a\circ_\lambda b&=& a\circ b,\\
  {[a_\lambda b]}&=& \partial(b\circ D (a))+\lambda(a\circ D (b)+b\circ D(a))+[b,a],\quad\forall~a,b\in A.
  \end{eqnarray}
\end{pro}
\begin{proof}
Since $D$ is a derivation on the commutative algebra $(A,\circ)$, $(A,\ast)$  is a Novikov algebra (\cite{Gel}). Since $D$ is a derivation on the Lie algebra $(A,[\cdot,\cdot])$, by the Leibniz rule of the Poisson algebra, \eqref{eq:GD condition} follows. Thus $(A,\ast,[\cdot,\cdot])$ is a $\GD$-algebra. By the commutative and associative rules of $(A,\circ)$, \eqref{eq:quad Poi1} follows. Also by the fact that $D$ is a derivation on the commutative algebra $(A,\circ)$, \eqref{eq:quad Poi2} follows. Thus $(A,\ast,\circ,[\cdot,\cdot])$ is a $\PGD$-algebra. By Proposition \ref{pro:quadratic comformal Poisson}, $(P,\circ_\lambda,[\cdot_\lambda\cdot])$ is a Poisson conformal algebra.
\end{proof}

\begin{cor}\label{cor:CAD}
 Let $(A,\circ)$ be a commutative associative algebra with a derivation $D$.  Then $(A,\ast,\circ,[\cdot,\cdot]=0)$ is a $\PGD$-algebra,  where the operation $\ast$ is given by \eqref{eq:Novikov product}. Furthermore, $(P,\circ_\lambda,[\cdot_\lambda\cdot])$ is a Poisson conformal algebra, where $\circ_\lambda$ and $[\cdot_\lambda\cdot]$ are given by
  \begin{eqnarray}
  a\circ_\lambda b&=& a\circ b,\\
  {[a_\lambda b]}&=& \partial(b\circ D (a))+\lambda(a\circ D (b)+b\circ D(a)),\quad\forall~a,b\in A.
  \end{eqnarray}
\end{cor}

Let $A=\Comp[x_1,x_2,\cdots,x_n]$ be the algebra of polynomials in $n$ variables. Let $$\huaD_n=\{\partial_{x_1},\partial_{x_2},\ldots,\partial_{x_n}\}$$ be the system of derivations over $A$. For any  polynomial $f\in A$, the endomorphisms
  $$f\partial_{x_i}:A\longrightarrow A,\quad (f\partial_{x_i})(g)=f\partial_{x_i}(g),\quad \forall~g\in A$$
  are derivations of $A$. Denote by $A\huaD_n=\{\sum_{i=1}^nf_i\partial_{x_i}\mid f_i\in A,\partial_{x_i}\in\huaD_n\}$ the space of derivations.

\begin{ex}
Let $A=\Comp[x_1,x_2,\cdots,x_n]$ be the algebra of polynomials in $n$ variables. For any $D\in A\huaD_n$, define $\circ:A\times A\longrightarrow A$ and  $\ast:A\times A\longrightarrow A$ by
  \begin{eqnarray*}
  g\circ h&=& gh,\\
g\ast h&=&gD(h)=\sum_{i=1}^n gf_i\partial_{x_i}(h),\quad\forall~g,h\in A.
  \end{eqnarray*}
  Then by Proposition \ref{cor:CAD}, $(A,\ast,\circ,[\cdot,\cdot]=0)$ is a $\PGD$-algebra. Furthermore, $(P,\circ_\lambda,[\cdot_\lambda\cdot])$ is a Poisson conformal algebra, where $\circ_\lambda$ and $[\cdot_\lambda\cdot]$ are given by
  \begin{eqnarray}
  g\circ_\lambda h&=& gh,\\
  {[g_\lambda h]}&=& \sum_{i=1}^n\partial( gf_i\partial_{x_i}(h))+\sum_{i=1}^n\lambda( gf_i\partial_{x_i}(h)+hf_i\partial_{x_i}(g)),\quad\forall~g,h\in A.
  \end{eqnarray}
\end{ex}

\begin{ex}\label{ex:comforal2}
	Let  $A=\Comp[x]$ be the algebra of polynomials in one variable. Define $\circ:A\times A\longrightarrow A$ and  $\ast:A\times A\longrightarrow A$ by
	\begin{eqnarray*}
		x^m\circ x^n&=&x^{m+n},\\
		x^m\ast x^n&=&x^m\partial_x(x^n)= nx^{m+n-1}.
	\end{eqnarray*}
	Then  $(A,\ast,\circ,[\cdot,\cdot]=0)$ is a $\PGD$-algebra. Furthermore, $(A,\circ_\lambda,[\cdot_\lambda\cdot])$ is a Poisson conformal algebra, where $\circ_\lambda$ and $[\cdot_\lambda\cdot]$ are given by
	\begin{eqnarray*}
		x^m\circ_\lambda x^n= x^{m+n},\quad	{[x^m_\lambda x^n]}= (n\partial +(n+m)\lambda)x^{n+m-1}.
	\end{eqnarray*}
\end{ex}

\section{Formal distribution (noncommutative)  Poisson algebras and (noncommutative)  Poisson conformal algebras}\label{sec: distribution}

Let $ V $ be a vector space. A {\bf $ V $-valued formal distribution} is a series of the form $ a(z)=\sum_{n\in\mathbb{Z}}a_{n}z^{-n-1} $, where $ a_{n}\in V $ and $ z $ is an indeterminate. We denote the space of such distributions by $ V[[z^{\pm1}]] $. Vectors $a_n$ is called {\bf Fourier coefficients} of $a(z)$.

Let $(A,[\cdot,\cdot])$ be a Lie  algebra.  Let $a$ and $b$ be $A$-valued formal distributions. Define the bracket of $a(z)$ and $b(w)$ by
  $$ [a(z), b(w)]=\sum_{m,n\in\mathbb{Z}}[ a_{m},b_{n}]z^{-m-1} w^{-n-1} .$$
We say the pair $(a,b)$ is a {\bf local pair} if for some $N\in\Nat$, the following relation holds:
\begin{eqnarray*}
	(z-w)^{N}[a(z),b(w)]=0.
	\end{eqnarray*}
Let $(A,[\cdot,\cdot])$ be a Lie  algebra and $ \huaF $  a family of pairwise local $ A $-valued formal distributions.  If the Fourier coefficients of formal distributions from $ \huaF $ span $ A $, then the pair $ ((A,\circ),\huaF) $ is called a {\bf formal distribution Lie algebra}.

Let $(A,\circ)$ be an (commutative) associative algebra.  Let $a$ and $b$ be $A$-valued formal distributions. Define the multiplication of $a(z)$ and $b(w)$ by
  $$ a(z)\circ b(w)=\sum_{m,n\in\mathbb{Z}}(a_{m}\circ b_{n})z^{-m-1}w^{-n-1}.$$
We say $(a,b)$ is a {\bf local pair} if for some $N\in\Nat$, the following relations hold:
\begin{eqnarray*}
	(z-w)^{N}a(z)\circ b(w)=0 ,\quad (z-w)^{N}b(z)\circ a(w)=0.
	\end{eqnarray*}
Let $(A,\circ)$ be an (commutative) associative  algebra and $ \huaF $ be a family of pairwise local $ A $-valued formal distributions.  If the Fourier coefficients of formal distributions from $ \huaF $ span $ A $, then the pair $ ((A,\circ),\huaF) $ is called a {\bf formal distribution (commutative) associative algebra}.
	
\begin{defi}
Let $(P,\circ,[\cdot,\cdot]) $ be a (noncommutative) Poisson algebra.
 \begin{itemize}
  \item[{\rm(i)}]   Let $a$ and $b$ be $P$-valued formal distributions. A pair $(a,b)$  is called {\bf local} over the (noncommutative) Poisson algebra $P$ if  $(a,b)$ is both a local pair on the Lie algebra $(P,[\cdot,\cdot])$ and a local pair on the  associative algebra $(P,\circ)$.
  \item[{\rm(ii)}] Let $ \huaF $ be  a family of pairwise local $ P $-valued formal distributions.  If $ ((P,[\cdot,\cdot]),\huaF) $ is a formal distribution Lie algebra and $ ((P,\circ),\huaF) $ is a formal distribution (commutative) associative algebra,
  then the pair $ ((P,\circ,[\cdot,\cdot]),\huaF) $ is called a {\bf formal distribution (noncommutative) Poisson algebra}.
      \end{itemize}
\end{defi}

For $ a(z)=\sum_{n\in\mathbb{Z}}a_{n}z^{-n-1}\in U[[z^{\pm1}]] $, we can define
\begin{eqnarray*}
	{\rm Res}_{z}a(z)=a_{0}.
\end{eqnarray*}
Let $(a,b)$ be a local pair over the Lie algebra (or (commutative) associative algebra).  There exists $N\in\Nat$ such that $(z-w)^{N}[a(z), b(w)]=0~((z-w)^{N}(a(z)\circ b(w))=0)$. By Corollary 2.2 in \cite{Kac}, we have
\begin{eqnarray*}
	[ a(z), b(w)]=\sum_{n=0} ^{N-1}(a_{[n]}b)(w)\partial^{(n)}_w\delta(z-w),\quad 
		a(z)\circ b(w)=\sum_{n=0} ^{N-1}(a_{(n)}b)(w)\partial^{(n)}_w\delta(z-w),
\end{eqnarray*}
where $a_{[n]}b={\rm Res}_{z}(z-w)^{n}[a(z), b(w)]$, $a_{(n)}b={\rm Res}_{z}(z-w)^{n}a(z)\circ b(w)$,  $\partial_w^{(n)}=\frac{\partial^n_w}{n!}$ and $\delta(z-w)=\sum_{n\in\mathbb{Z}}z^{-n-1}w^{n}$.

 The formal distribution $(a_{[n]}b)$ ($a_{(n)}b$) is $\Comp$-bilinear in $a$ and $b$ and it is called the $n$-th product of $a$ and $b$.

 \emptycomment{By a direct calculation, we have
 \begin{eqnarray*}
 (\partial a)_{[n]} b=-na_{[n-1]}b,\quad a_{[n]}\partial b=\partial(a_{[n]}b)+na_{[n-1]}b,\\
 	(\partial a)_{(n)} b=-na_{(n-1)}b,\quad a_{(n)}\partial b=\partial(a_{(n)}b)+na_{(n-1)}b.
 \end{eqnarray*} }

 Let $(A,\huaF)$ be a formal distribution Lie algebra or (commutative) associative algebra. Denote by $ \bar{\huaF} \subseteq A[[z^{\pm1}]]$  a minimal $ \mathbb{C}[\partial_z]$-module  containing $ \huaF $ which is closed under all $ j $-th products for $j\in\Nat$.

 Next, we recall Dong's lemma for Lie algebras and (commutative) associative algebras.

 \begin{lem}{\rm(Dong's lemma)} Let $a$, $b$ and $c$ be pairwise mutually local formal distributions over  a Lie algebra or (commutative) associative algebra. Then for all $n\in\Nat$, the formal distributions $a_{[n]}b~(a_{(n)}b)$ and $c$ are mutually local. 	
 \end{lem}
 \begin{proof}
 See the proof of Lemma 2.8 in \cite{Kac}.	
 \end{proof}

 By Dong's lemma, if $(A,\huaF)$ is a formal distribution Lie algebra or (commutative) associative algebra, then $\bar{\huaF}$ consists of mutually local formal distributions. Note that $\bar{\huaF}$ is a $ \mathbb{C}[\partial]$-module, where the the action of $\partial$ is given by $(\partial(a))(z)=\partial_z a(z)$.

 \begin{thm}{\rm(\cite{Kac97})}\label{thm:Lie-Ass conformal}
  \begin{itemize}
  \item[{\rm(i)}] Let $ ((A,[\cdot,\cdot]),\huaF) $ be a formal distribution Lie algebra. Define a $ \lambda $-bracket over $ \bar{\huaF} $
 	\begin{eqnarray}\label{eq:Lie con}
 		[a_{\lambda}b]=\sum_{n\in\Nat}\lambda^{(n)}a_{[n]}b\in\bar{\huaF}[\lambda],\quad \forall ~a,b\in\bar{\huaF}.
 	\end{eqnarray}
  Then $ \bar{\huaF} $ becomes a Lie conformal algebra under this $ \lambda $-bracket.
  \item[{\rm(ii)}] Let $ ((A,\circ),\huaF) $ be a formal distribution (commutative) associative algebra. Define a $ \lambda $-multiplication over $ \bar{\huaF} $
 	\begin{eqnarray}\label{eq:ass con}
 		a\circ_{\lambda}b=\sum_{n\in\Nat}\lambda^{(n)}a_{(n)}b\in\bar{\huaF}[\lambda],\quad \forall ~a,b\in\bar{\huaF}.
 	\end{eqnarray}
 Then $ \bar{\huaF} $ becomes an (commutative) associative conformal algebra under this $ \lambda $-multiplication.
    \end{itemize}
 \end{thm}

Recall that the {\bf formal Fourier transform} of a formal distribution $ a(z,w)\in U[[z^{\pm1},w^{\pm1}]] $ by the formula:
\begin{eqnarray*}
	F^{\lambda}_{z,w}(a(z,w))={\rm Res}_{z}e^{\lambda(z-w)}a(z,w),
\end{eqnarray*}
where $ e^{\lambda(z-w)}=\sum_{j\in\Nat}\lambda^{(j)}(z-w)^{j}\in \mathbb{C}[[z^{\pm1},w^{\pm1}]][[\lambda]] $.

\begin{lem}{\rm(\cite{Kac})}\label{lem:Fourier}
	 For $ a(z,w)\in U[[z^{\pm1},w^{\pm1}]] $ or $ a(z,w,x)\in U[[z^{\pm1},w^{\pm1},x^{\pm1}]] $, then we have
	\begin{itemize}
	\item [{\rm(i)}]  $ F^{\lambda}_{z,w}\partial_{z}(a(z,w))=-\lambda F^{\lambda}_{z,w}(a(z,w))=[\partial_{w},F^{\lambda}_{z,w}](a(z,w)) $,
	\item [{\rm(ii)}]  $ F^{\lambda}_{z,w}F^{\mu}_{x,w}a(z,w,x)=F^{\lambda+\mu}_{x,w}F^{\lambda}_{z,x}a(z,w,x) $.
	\end{itemize}
\end{lem}

\begin{thm}
	Suppose that $ (P,\huaF) $ is a formal distribution (noncommutative) Poisson algebra. 	Then $ (\bar{\huaF},\circ_\lambda,[\cdot_\lambda \cdot]) $ is a (noncommutative) Poisson conformal algebra, where the $\lambda$-operation  $\circ_\lambda$ is given by \eqref {eq:ass con}  and the $\lambda$-bracket $[\cdot_\lambda \cdot]$ is given by  \eqref {eq:Lie con}. We denote this (noncommutative)Poisson conformal algebra by ${\rm Con}(P,\huaF)$.
\end{thm}
\begin{proof}Since $ (P,\huaF) $ is a formal distribution (noncommutative)  Poisson algebra, by Theorem \ref{thm:Lie-Ass conformal}, $ (\bar{\huaF},\circ_\lambda) $ is an (commutative) associative algebra and $ (\bar{\huaF},[\cdot_\lambda \cdot]) $ is a Lie conformal algebra. Furthermore, it was shown in \cite{Kac} that 	the (commutative) associative conformal $\lambda$-operation  $\circ_\lambda$ and  the Lie conformal  $\lambda$-bracket $[\cdot_\lambda \cdot]$ can be written by
$$ a(w)\circ_{\lambda}b(w)=F_{z,w}^{\lambda}a(z)b(w) ,\quad  [a(w)_{\lambda}b(w)]=F_{z,w}^{\lambda}[a(z),b(w)]. $$
In the following, we show that the conformal Leibniz rule holds. By definition of formal Fourier transformations and  Lemma \ref{lem:Fourier}, we have
\begin{eqnarray*}
[a(w)_{\lambda}(b(w)\circ_{\mu}c(w))]&=&[a(w)_{\lambda}(F_{x,w}^{\mu}(b(x)\circ c(w))]=F_{z,w}^{\lambda}[a(z),F_{x,w}^{\mu}(b(x)\circ c(w))]\\
	&=&F_{z,w}^{\lambda}F_{x,w}^{\mu}[a(z),b(x)\circ c(w)];\\
	{[a(w)_{\lambda}b(w)]}\circ_{\lambda+\mu} c(w)&=&F_{x,w}^{\lambda+\mu}[a(x)_{\lambda}b(x)]\circ c(w)=F_{x,w}^{\lambda+\mu}F_{z,x}^{\lambda}[a(z),b(x)]\circ c(w)\\
	&=&F_{z,w}^{\lambda}F_{x,w}^{\mu}[a(z),b(x)]\circ c(w);\\
	b(w)\circ_{\mu}[a(w)_{\lambda}c(w)]&=&F_{x,w}^{\mu}b(x)\circ [a(w)_{\lambda}c(w)]=F_{x,w}^{\mu}b(x)\circ (F_{z,w}^{\lambda}[a(z),c(w)])\\
&=&F_{z,w}^{\lambda}F_{x,w}^{\mu}b(x)\circ [a(z),c(w)].
\end{eqnarray*}
Thus we have
\begin{eqnarray*}
&&[a(w)_{\lambda}(b(w)_{\mu}\circ c(w))]-[a(w)_{\lambda}b(w)]\circ_{\lambda+\mu}c(w)-b(w)\circ_{\mu}[a(w)_{\lambda}c(w)]\\
&=&F_{z,w}^{\lambda}F_{x,w}^{\mu}\big([a(z),b(x)\circ c(w)]-[a(z),b(x)]\circ c(w)-b(x)\circ [a(z),c(w)] \big)\\
&=&0.
\end{eqnarray*}
The last equality follows form Leibniz rule of the (noncommutative) Poisson algebra $P$.

Therefore, $ (\bar{\huaF},\circ_\lambda,[\cdot_\lambda \cdot]) $ is a (noncommutative) Poisson conformal algebra.
 \end{proof}

Let $A$ be a vector space.  Denote  {\rm Coeff} $A$  by the quotient of the vector space with the basis $ a_{n}~(a\in A, n\in\mathbb{Z})$ by the subspace spanned over $ \mathbb{C} $ by elements:
	\begin{eqnarray*}
		(\alpha a)_{n}-\alpha a_{n},\quad (a+b)_{n}-a_{n}-b_{n},\quad (\partial a)_{n}+na_{n-1},\quad  \forall ~ a,b\in A, \alpha\in\mathbb{C}, n\in\mathbb{Z}.
	\end{eqnarray*}
	
Furthermore, if $A$ is  a Lie conformal algebra, then $A$ equipped with the following bracket
\begin{eqnarray}
    	\label{eq:j1}[a_{m},b_{n}]&=&\sum_{j\in\Nat}\binom{m}{j}(a_{[j]}b)_{m+n-j}
\end{eqnarray}
is a Lie algebra on {\rm Coeff} $A$ , which is called a {\bf coefficient Lie algebra}.

If $A$ is  an (commutative) associative conformal algebra, then $A$ equipped with the following operation
\begin{eqnarray}
    	\label{eq:j2}a_{m}\circ b_{n}&=&\sum_{j\in\Nat}\binom{m}{j}(a_{(j)}b)_{m+n-j}
    	\end{eqnarray}
is an (commutative) associative algebra on {\rm Coeff} $A$, which is called a {\bf coefficient (commutative) associative algebra}.	

See \cite{Kac97,Kac} for more details on coefficient Lie algebras and coefficient (commutative) associative algebras.

\begin{pro}\label{pro:Pa}
	Let $ (P,\circ_\lambda,[\cdot_\lambda\cdot]) $ be a (noncommutative) Poisson conformal algebra. Then $({\rm Coeff}~P,\circ,[\cdot,\cdot])$ is a (noncommutaive) Poisson algebra, where
	the operation $\circ$ is given by \eqref{eq:j2} and the bracket $[\cdot,\cdot]$ is given by \eqref{eq:j1}. We call this (noncommutaive) Poisson algebra a {\bf coefficient (noncommutative) Poisson algebra}.
	\end{pro}
	\begin{proof}
	Since 	$(P,\circ_\lambda)$ is an (commutative) associative algebra, $({\rm Coeff}~P,\circ)$ is an (commutative) associative algebra. Similarly, $({\rm Coeff}~P,[\cdot,\cdot])$ is a Lie algebra. We only need to show that the  Leibniz rule holds. By the conformal Leibniz rule, we have
	\begin{eqnarray*}
&&[a_{m},b_{n}\circ c_{k}]-b_{n}\circ[a_{m},c_{k}]-[a_{m},b_{n}]\circ c_{k}\\
&=&\sum_{i,j\in\Nat}\binom{m}{i}\binom{n}{j}(a_{[i]}(b_{(j)}c))_{m+n+k-i-j}-\sum_{i,j\in\Nat}\binom{m}{i}\binom{n}{j}(b_{(j)}(a_{[i]}c))_{m+n+k-i-j}\\
&&-\sum_{i',j'\in\Nat}\binom{m}{i'}\binom{m+n-i'}{j'}((a_{[i']}b)_{(j')}c)_{m+n+k-i'-j'}\\
&=&\sum_{i,j\in\Nat}\binom{m}{i}\binom{n}{j}\Big(\sum_{s=0}^{i}\binom{i}{s}(a_{[s]}b)_{(i+j-s)}c+b_{(j)}(a_{[i]}c)\Big)_{m+n+k-i-j}-\sum_{i,j\in\Nat}\binom{m}{i}\binom{n}{j}(b_{(j)}(a_{[i]}c))_{m+n+k-i-j}\\
&&-\sum_{i',j'\in\Nat}\binom{m}{i'}\binom{m+n-i'}{j'}((a_{[i']}b)_{(j')}c)_{m+n+k-i'-j'}\\
&=& \sum_{i,j\in\Nat}\binom{m}{i}\binom{n}{j}\sum_{s=0}^{i}\binom{i}{s}((a_{[s]}b)_{(i+j-s)}c)_{m+n+k-i-j}-\sum_{i',j'\in\Nat}\binom{m}{i'}\binom{m+n-i'}{j'}((a_{[i']}b)_{(j')}c)_{m+n+k-i'-j'}\\
&=&\sum_{i,i',j'\in\Nat}\binom{m}{i}\binom{n}{i'+j'-i}\binom{i}{i'}((a_{[i']}b)_{(j')}c)_{m+n+k-i'-j'}-\sum_{i',j'\in\Nat}\binom{m}{i'}\binom{m+n-i'}{j'}((a_{[i']}b)_{(j')}c)_{m+n+k-i'-j'}\\
&=&0.
	\end{eqnarray*}
	The last equality follows from
	\begin{eqnarray*}
	\sum_{i,i',j'\in\Nat}\binom{m}{i}\binom{n}{i'+j'-i}\binom{i}{i'}
	&=&\sum_{i,i',j'\in\Nat}\binom{m}{i'}\binom{m-i'}{i-i'}\binom{n}{i'+j'-i}\\
	&=&\sum_{s,i',j'\in\Nat}\binom{m}{i'}\binom{m-i'}{s}\binom{n}{j'-s}\\
	&=&\sum_{i',j'\in\Nat}\binom{m}{i'}\binom{m+n-i'}{j'}.
\end{eqnarray*}
Thus $({\rm Coeff}~P,\circ,[\cdot,\cdot])$ is a (noncommutaive) Poisson algebra.
	\end{proof}

\begin{pro}
Let $ (P,\circ_\lambda,[\cdot_\lambda\cdot]) $ be a (noncommutative) Poisson conformal algebra. Define $D:P\rightarrow P$ by
$$D(a_n)=-na_{n-1}.$$	
Then $D$ is a derivation on the (noncommutaive) Poisson algebra $({\rm Coeff}~P,\circ,[\cdot,\cdot])$.
\end{pro}	
	\begin{proof}
		It follows by a direct calculation.
	\end{proof}	

Furthermore, by Proposition \ref{ex:derivation-HMA}, we have	
\begin{ex}
Let $ (P,\circ_\lambda,[\cdot_\lambda\cdot]) $ be a  Poisson conformal algebra. Then there is a new Poisson conformal algebra $ ({\rm Coeff}~P,\tilde{\circ}_\lambda,\{\cdot_\lambda\cdot\}) $, where $\tilde{\circ}_\lambda$ and $\{\cdot_\lambda\cdot\}$ are given by
  \begin{eqnarray*}
	a_m\tilde{\circ}_\lambda b_n&=& \sum_{j\in\Nat}\binom{m}{j}(a_{(j)}b)_{m+n-j},\\
	\{(a_m)_\lambda (b_n)\}&=& -\sum_{j\in\Nat}\Big(m\binom{n}{j}  \partial( (b_{(j)}a)_{m+n-j-1}   ) +\lambda \big( n\binom{m}{j}   (a_{(j)}b)_{m+n-j-1} \\
	&& +  m\binom{n}{j} (b_{(j)}a)_{m+n-j-1} \big)       \Big)+\sum_{j\in\Nat}\binom{m}{j}(b_{[j]}a)_{m+n-j}.
\end{eqnarray*}
	\end{ex}

\begin{thm}
	Suppose that $ P $ is a (noncommutative) Poisson conformal algebra. Let $ \tilde{P}=P[t^{\pm1}]/(\partial+\partial_{t})P[t^{\pm1}] $ and  $ a_{n} $ denote the image of $ at^{n} $ in $  \tilde{P}$.
	\begin{itemize}
    \item [{\rm(i)}] Define the bracket $[\cdot,\cdot]$ and the operation $\circ$ on $  \tilde{P}$ by    \begin{eqnarray*}
    	[a_{m},b_{n}]=\sum_{j\in\Nat}\binom{m}{j}(a_{[j]}b)_{m+n-j},\quad a_{m}\circ b_{n}=\sum_{j\in\Nat}\binom{m}{j}(a_{(j)}b)_{m+n-j}.
    \end{eqnarray*}
Then $( \tilde{P},\circ,[\cdot,\cdot]) $ is a  (noncommutative) Poisson algebra.
    \item [{\rm(ii)}] Set $ \huaF=\{a(z)=\sum_{n\in\mathbb{Z}}a_{n}z^{-n-1}\arrowvert a\in P\}\subseteq \tilde{P}[[z^{\pm1}]] $. Then $ \huaF $ is a family of pairwise local $ \tilde{P} $-valued formal distributions and thus the pair $ (\tilde{P},\huaF) $ is a formal distribution (noncommutative) Poisson algebra.
    \item [{\rm(iii)}] The formal distribution  (noncommutative) Poisson conformal algebra $ {\rm Con}~(\tilde{P},\huaF)  $ is isomorphism to the (noncommutative) Poisson conformal algebra $P$.
	\end{itemize}
\end{thm}
\begin{proof}
	(i) Note that $(\partial+\partial_{t})P[t^{\pm1}]={\rm span}_{\Comp}\{(\partial a)t^n+nat^{n-1}\mid a\in P,n\in\Integ\}$.  It is not hard to check that
	\begin{eqnarray*}
		(\alpha a)_{n}=\alpha a_{n},\quad (a+b)_{n}=a_{n}+b_{n},\quad (\partial a)_{n}=-na_{n-1},\quad  \forall ~ a,b\in A, \alpha\in\mathbb{C}, n\in\mathbb{Z}.
	\end{eqnarray*}
	Thus $\tilde{P}$ is just {\rm Coeff} $P$. By Proposition \ref{pro:Pa}, $( \tilde{P},\circ,[\cdot,\cdot] $ is a  (noncommutative) Poisson algebra.
	
	(ii) It was shown in \cite{Kac} that $ \huaF $ is a family of pairwise local $ \tilde{P} $-valued formal distributions over the both Lie algebra and (commutative) associative algebra. Then 	$ \huaF $ is a family of pairwise local $ \tilde{P} $-valued formal distributions	over the  (noncommutative) Poisson algebra. Thus the pair $ (\tilde{P},\huaF) $ is a formal distribution (noncommutative) Poisson algebra.
	
	(iii) Since the map $P\rightarrow\tilde{P}:a\mapsto a(z)$ is both a Lie conformal algebra isomorphism and an associative conformal algebra isomorphism,  the formal distribution  (noncommutative) Poisson conformal algebra $ {\rm Con}~(\tilde{P},\huaF)  $ and the (noncommutative) Poisson conformal algebra $P$ are isomorphism.
	\end{proof}
	
	\begin{ex}
		The coefficient Poisson algebra $({\rm Coeff}~A,\circ,[\cdot,\cdot])$ corresponding to the Poisson conformal algebra in Proposition \ref{ex:derivation-HMA} is given by
		\begin{eqnarray*}
			a_m\circ b_n&=&(a\circ b)_{m+n},\\
			{[a_m,b_n]}&=&([a,b])_{m+n}+(ma\circ D(b)-nb\circ D(a))_{m+n-1},\quad \forall~a,b\in A.
		\end{eqnarray*}
	\end{ex}
	
	\begin{ex}
		The coefficient Poisson algebra $({\rm Coeff}~A,\circ,[\cdot,\cdot])$ corresponding to the Poisson conformal algebra in Example \ref{ex:comforal2} is given by
		\begin{eqnarray*}
			(x^k)_m\circ (x^l)_n&=&(x^{k+l})_{m+n},\\
			{[(x^k)_m,(x^l)_n]}&=&(lm-kn)(x^{k+l-1})_{m+n-1}.
		\end{eqnarray*}
	\end{ex}

\section{Poisson conformal algebras and deformation quantization via  associative conformal algebras}\label{sec:deformation quantization}
In this section, we introduce the notion of conformal formal deformations of
commutative associative conformal algebras  and show that Poisson conformal algebras  are the corresponding semi-classical limits.  Furthermore,  we show that conformal infinitesimal deformation and extension of conformal $n$-deformation to conformal $(n+1)$-deformation of a commutative associative conformal algebra are classified by the second and the third cohomology groups of the associative conformal algebra respectively.

\begin{defi}
	A {\bf module} $ V $ over an associative conformal algebra $ A $ is a $ \mathbb{C}[\partial] $-module endowed with a pair $ \mathbb{C} $-bilinear maps $ A\times V \rightarrow V[\lambda]: ~(a,v) \mapsto a\circ_{\lambda}v $ and  $ V\times A \rightarrow V[\lambda]:~(v,a) \mapsto v\circ_{\lambda} a $, satisfying
	\begin{eqnarray*}
		(\partial a)\circ_\lambda v&=&-\lambda (a\circ_\lambda v),\quad a\circ_\lambda (\partial v)=(\partial+\lambda)(a\circ_\lambda v),\\
		(\partial v)\circ_\lambda a&=&-\lambda (v \circ_\lambda a),\quad v \circ_\lambda(\partial a)=(\partial+\lambda)(v\circ_\lambda a),\\
		(a\circ_\lambda b)\circ_{\lambda +\mu}v&=&a\circ_\lambda(b\circ_\mu v),\quad (v\circ_\lambda b)_{\lambda+\mu}a=v\circ_\lambda(b\circ_\mu a),\\
		(a\circ_\lambda v)\circ_{\lambda+\mu}b&=&a\circ_\lambda(v\circ_\mu b),\quad \forall~a,b\in A,v\in V.
	\end{eqnarray*}
\end{defi}

For convenience, we set $l(a)_\lambda v=a\circ_\lambda v$ and $r(a)_\lambda v=v\circ_{-\lambda-\partial}a$. Then a structure of a module $V$ over an associative conformal algebra $A$ can be equivalently described by two $\Comp[\partial]$-module homomorphisms $l$ and $r$ from $A$ to ${\rm gc}(V)$. Here ${\rm gc}(V)$ represents the set of conformal linear transformations from  $V$ to $V$. More precisely,
$${\rm gc}(V):=\{ a: V\rightarrow V[\lambda] | ~ a_\lambda(\partial v)=(\partial+\lambda)a_\lambda v,~\forall~v\in V  \}.$$

\begin{ex}
Let $(V;l,r)$ be a finite module over an associative conformal algebra $A$. Let $l^*,r^*$ be two $\Comp[\partial]$-module homomorphisms  from $A$ to $V^{c*}$ defined by
$$(l(a)_\lambda\varphi)_\mu v=\varphi_{\mu-\lambda}(l(a)_\lambda v),~ (r(a)_\lambda\varphi)_\mu v=\varphi_{\mu-\lambda}(r(a)_\lambda v) $$
for $a\in A,\varphi\in V^{c*}$ and $v\in V$. Then $(V^{*c};r^*,l^*)$ is a  module over $A$.
	\end{ex}

\begin{ex}
Let $A$ be an associative conformal algebra. Define two $\Comp[\partial]$-module homomorphisms $L$ and $ R$ from $A$ to ${\rm gc}(V)$ by $L(a)_\lambda b=a_\lambda b$ and $R(a)_\lambda b=b_{-\lambda-\partial}a$ for $a,b\in A$. Then $(A;L,R)$ is a module over $A$. Furthermore, if $A$ is a finite associative conformal algebra, then $(A^{*c};R^*,L^*)$ is also a module over $A$.

\end{ex}

\emptycomment{\begin{lem}
		Suppose that $ A $ is an associative conformal algebra and $ V $ is an $ A $-module. For all $ a, b\in A, ~v\in Vv $, we have
		\begin{eqnarray}
			\label{AM1}a\circ_{\lambda}(v_{-\partial-\mu}\circ b)&=&(a\circ_{\lambda}v)\circ_{-\partial-\mu} b,\\
			\label{AM2}v\circ_{-\partial-\lambda} (a\circ_{\mu}b)&=&(v\circ_{-\partial-\mu} a)\circ_{-\partial+\mu-\lambda} b.
		\end{eqnarray}
	\end{lem}
	\begin{proof}
		It follows by a direct calculations.
\end{proof}}

The following conclusion is direct.
\begin{pro}
	Let $A$ be an associative conformal algebra and let $V$ be a $\Comp [\partial]$-module endowed with a pair $ \mathbb{C} $-bilinear maps $ A\times V \rightarrow V[\lambda]: ~(a,v) \mapsto a\circ_{\lambda}v $ and  $V\times A \rightarrow V[\lambda]:~(v,a) \mapsto v\circ_{\lambda} a $. Then $V$ is a module over the associative conformal algebra $A$ if and only if $A\oplus V$
	endowed with a $\Comp[\partial]$-module structure and a $\lambda$-operation defined by:
	\begin{eqnarray*}
		\partial(a+v)&=&\partial a+\partial v, \\
		{(a+u)\circ_{\lambda}(b+v)}&=&a\circ_{\lambda}b+a\circ_{\lambda}v+u\circ_{\lambda}b ,\quad \forall~a,b\in A,u,v\in V
	\end{eqnarray*}
	is an associative conformal algebra. This associative conformal algebra is called the {\bf semi-direct product} of $A$ and $V$.
\end{pro}

Let $A$ and $V$ be $\Comp[\partial]$-modules. For $n\geq 1 $, an {\bf $n$-$\lambda$-bracket} on $A$ with coefficients in $V$ is a $\Comp$-linear map $\gamma:A^{\otimes n}\rightarrow \Comp[\lambda_{1},\cdots,\lambda_{n-1}]\otimes V$, denoted by
$$a_{1}\otimes\cdots\otimes a_{n}\mapsto \{a_{1\lambda_{1}}\cdots a_{n-1 \lambda_{n-1}}a_{n}\}_{\gamma}$$
satisfying the following conditions:
\begin{eqnarray}
	\label{eq:deformCA}  \{a_{1\lambda_{1}}\cdots(\partial a_{i})_{\lambda_{i}}\cdots a_{n-1 \lambda_{n-1}}a_{n}\}_{\gamma}&=&-\lambda_{i}\{a_{1\lambda_{1}}\cdots a_{n-1 \lambda_{n-1}}a_{n}\}_{\gamma},\quad 1\leq i\leq n-1,\\
	\label{eq:deforLA}\{a_{1\lambda_{1}}\cdots a_{n-1 \lambda_{n-1}}(\partial a_{n})\}_{\gamma}&=&(\lambda_{1}+\cdots +\lambda_{n-1}+\partial)\{a_{1\lambda_{1}}\cdots a_{n-1\lambda_{n-1}}a_{n}\}_{\gamma}.
\end{eqnarray}

Let $A$  be an  associative conformal algebra and $V$ an $A$-module. For $n\geq1$, we denote by $C_{\rm H}^{n}(A,V)$ the space of all $n$-$\lambda$-brackets on $A$ with coefficients in $V$. Set $C_{\rm H}^{\bullet}(A,V)=\oplus _{n\in\Nat}C_{\rm H}^{n}(A,V)$.

\emptycomment{For $v+\partial V \in C_{\rm H}^{0}(A,V)$, we let $\dM_{\rm H} (v+\partial V)\in C_{\rm H}^{1}(A,V)$ be the following $\Comp[\partial]$-module homomorphism
\begin{equation}
	\{a\}_{\dM_{\rm H} (v+\partial V)}: =a\circ_{-\partial}v-v\circ_{-\partial} a,\quad \forall~a\in A.
\end{equation}}
For any $\gamma\in C_{\rm H}^{n}(A,V)$ with $n\geq 1$, the corresponding conformal Hochschild differential $\dM_{\rm H} \gamma \in C_{\rm H}^{n+1}(A,V)$ is defined by
\begin{eqnarray}\label{eq:d}
	\nonumber\{{a_{1}}_{\lambda_{1}}\dots {a_{n}}_{\lambda_{n}}{ a_{n+1}}\}_{\dM_{\rm H}\gamma}
	&:=&{a_{1}}\circ_{\lambda_{1}}\{{a_{2}}_{\lambda_{2}}\dots {a_{n}}_{\lambda_{n}}{ a_{n+1}}\}_{\gamma}\\
	&&+\sum_{i=1}^{n}(-1)^{i} \{{a_{1}}_{\lambda_{1}}\dots {a_{i-1}}_{\lambda_{i-1}}({a_{i}}\circ_{\lambda_{i}} a_{i+1})_{\lambda_{i}+\lambda_{i+1}}{a_{i+2}}_{\lambda_{i+2}}\dots {a_{n}}_{\lambda_{n}}{a_{n+1}}\}_{\gamma}\\
	%\nonumber &&+(-1)^n\{{a_{1}}_{\lambda_{1}}\dots a_{n-1\lambda_{n-1}}({a_{n}}\circ_{\lambda_{n}}{a_{n+1}})\}_{\gamma}\\
	\nonumber&&+(-1)^{n+1}{\{{a_{1}}_{\lambda_{1}}\dots {a_{n-1}}_{\lambda_{n-1}}{ a_{n}}\}_{\gamma}}\circ_{\lambda^\flat_{n+1}} a_{n+1},
\end{eqnarray}
where $a_1,\cdots,a_{n+1}\in A$ and $\lambda_{n+1}^{\flat}=\sum_{j=1}^{n}\lambda_{j}$. We denote by $H_{\rm H}^{\bullet}(A,V)=\oplus _{k\in\Nat}H_{\rm H}^{k}(A,V)$ the cohomology of the associative conformal algebra $A$ with coefficients in $V$. See \cite{Dolg,Kol19} for more details on cohomology  of associative conformal algebras.

Let $A$ be an associaitve conformal algebra. It is well-known that $A[[\hbar]]$ is the completion of the $ ~\mathbb{\Comp}[\hbar] $-module with respect to the $\hbar$-adic topology. Formally we can extend the associative conformal algebra $\lambda$-operation $\circ_\lambda$ to a $ ~\mathbb{\Comp}[[\hbar]] $-bilinear  $\lambda$-operation  $ \circ^\hbar_{\lambda}$ on $ A[[\hbar]] $ by
\begin{equation}
	(\sum_{i\geq 0}a_i\hbar^i)\circ^\hbar_{\lambda}	(\sum_{j\geq 0}b_j\hbar^j)=\sum_{i,j\geq 0}(a_i\circ_\lambda b_j)\hbar^{i+j}\in A[[\hbar]] [[\lambda]].
\end{equation}
Because $	(\sum_{i\geq 0}a_i\hbar^i)\circ^\hbar_{\lambda}	(\sum_{j\geq 0}b_j\hbar^j)$ in general involves  infinite number of powers of $\lambda$, $(A[[\hbar]] ,\circ^\hbar_{\lambda}	)$ does not carry the structure of an  associative conformal algebra. However, note that for any $n\in \Nat$, the quotient space $A[[\hbar]] /\hbar^n A[[\hbar]]$ is naturally an associative conformal algebra. Motivated by this, we give the following notion of $\hbar$-adic associative conformal algebra.
\begin{defi}
	An {\bf $\hbar$-adic associative conformal algebra}  is a $ \mathbb{\Comp}[[\hbar]] $-module $A[[\hbar]] $, where $A$
	is a vector space, equipped with the $\hbar$-adic topology, a $ \mathbb{C}[\partial] $-module structure and a continuous $ \mathbb{\Comp}[[\hbar]] $-bilinear $\lambda$-operation  $ \circ^\hbar_{\lambda}:A[[\hbar]]\times A[[\hbar]] \rightarrow  A[[\hbar]] [[\lambda]]$  such that for any $n\in \Nat$, $(A[[\hbar]] /\hbar^n A[[\hbar]], \bar{\circ}^\hbar_{\lambda})$ is an associative conformal algebra over $ \mathbb{\Comp}[[\hbar]] $, where $\bar{\circ}^\hbar_{\lambda}$ is the natural quotient map of $\circ^\hbar_{\lambda}$, i.e., 
	$$(a+ \hbar^n A[[\hbar]])\bar{\circ}^\hbar_{\lambda}	(b+\hbar^n A[[\hbar]])=a\circ^\hbar_{\lambda} b+\hbar^n A[[\hbar]],\quad \forall~a,b\in A.$$
\end{defi}	
\begin{rmk}
	Let $(A[[\hbar]] , \circ^\hbar_{\lambda})$ be an $\hbar$-adic associative conformal algebra. For every pair of nonnegative integers $i$ and $j$ with $i\leq j$, the natural map $\pi_{ij}$ from $A[[\hbar]] /\hbar^j A[[\hbar]]$ onto $A[[\hbar]] /\hbar^i A[[\hbar]]$ is an associative conformal algebra homomorphism over $ \mathbb{\Comp}[[\hbar]] $. The associative conformal algebras $A[[\hbar]] /\hbar^n A[[\hbar]]$ together with these associative conformal algebra homomorphisms form an inverse system of associative conformal algebras over $ \mathbb{\Comp}[[\hbar]] $. Then the $\hbar$-adic associative conformal algebra $(A[[\hbar]] , \circ^\hbar_{\lambda})$  can be seen as an inverse limit of the associative conformal algebras $A[[\hbar]] /\hbar^n A[[\hbar]]$ over $ \mathbb{\Comp}[[\hbar]] $.
\end{rmk}

\begin{defi}
	Let $ A $ be a commutative associative conformal algebra.  A {\bf conformal  formal deformation} of $ A $ is 
	a continuous $ \mathbb{\Comp}[[\hbar]] $-bilinear $\lambda$-operation  $ \circ^\hbar_{\lambda}:A[[\hbar]]\times A[[\hbar]] \rightarrow  A[[\hbar]] [[\lambda]]$  such that $(A[[\hbar]] , \circ^\hbar_{\lambda})$ carries the structure of an $\hbar$-adic associative conformal algebra over $ \mathbb{\Comp}[[\hbar]] $, where $\circ^\hbar_{\lambda}$ extends naturally to a continuous $ \mathbb{\Comp}[[\hbar]] $-map on $A[[\hbar]]\otimes A[[\hbar]]$, such that for $a,b\in A$, 
		\begin{eqnarray*}
		a\circ^{\hbar}_{\lambda} b=\sum_{n=0}^{\infty}\hbar^{n}\{a_{\lambda}b\}_{\mu_{n}}=a\circ_{\lambda}b+\sum_{n=1}^{\infty}\hbar^{n}\{a_{\lambda}b\}_{\mu_{n}},
	\end{eqnarray*}
	where $ \mu_{n}: A\times A\rightarrow A[\lambda]~ $ for $ ~n\geqslant 0 $ are $\Comp$-bilinear maps with $ \mu_{0} $ being the commutative associative $\lambda$-operation on $ A $.
	\end{defi}

By the fact that $(A[[\hbar]] , \circ^\hbar_{\lambda})$ is an $\hbar$-adic associative conformal algebra, then we have
\begin{eqnarray}
\label{eq:semilinear}\{\partial a_\lambda b\}_{\mu_n}&=&-\lambda\{a_\lambda b\}_{\mu_n},\quad \{ a_\lambda \partial b\}_{\mu_n}~=~(\partial+\lambda)\lambda\{a_\lambda b\}_{\mu_n},\\
\label{fcd}\sum_{r+s=n}\{a_{\lambda}\{b_{\mu}c\}_{\mu_{s}}\}_{\mu_{r}}&=&\sum_{r+s=n}\{{\{a_{\lambda}b\}_{\mu_{s}}}_{\lambda+\mu}c\}_{\mu_{r}},\quad \forall ~n\geqslant 0.
\end{eqnarray}

\begin{thm}
	Let $ (A[[\hbar]],  \circ^\hbar_{\lambda}) $ be a conformal  formal deformation of $ (A,\circ_\lambda) $. Define
	\begin{eqnarray*}
		[a_{\lambda}b]= \{a_{\lambda}b\}_{\mu_{1}}-\{b_{-\lambda-\partial}a\}_{\mu_{1}},\quad \forall~a, b\in A.
	\end{eqnarray*}
	Then $ (A,\circ_\lambda,[\cdot_\lambda\cdot]) $ is a Poisson conformal algebra, which is called {\bf the semi-classical limit} of $ (A[[\hbar]], \circ^\hbar_{\lambda}) $. The $\hbar$-adic associative conformal algebra  $ (A[[\hbar]], \circ^\hbar_{\lambda}) $ is called a {\bf deformation quantization of $ (A,\circ_\lambda,[\cdot_\lambda\cdot]) $}.
	\end{thm}
	\begin{proof}
	Define the $\lambda$-bracket $ [\cdot_{\lambda}\cdot]_{\hbar} $ on $ A[[\hbar]] /\hbar^3 A[[\hbar]]$ by
\begin{eqnarray*}
	[a_{\lambda}b]_{\hbar}
	&=&a\bar{\circ}^\hbar_{\lambda}b-b\bar{\circ}^\hbar_{-\lambda-\partial}a,\quad \forall~a,b\in A.
\end{eqnarray*}
By the fact that $(A[[\hbar]] /\hbar^3 A[[\hbar]], \bar{\circ}^\hbar_{\lambda})$   is an associative conformal algebra, 
 $ (A[[\hbar]] /\hbar^3 A[[\hbar]],\bar{\circ}^\hbar_\lambda,[\cdot_\lambda\cdot]_\hbar) $ is a noncommutative Poisson conformal algebra. Thus, for $a,b,c\in A$, we have
\begin{eqnarray}
	\label{eq:cLr}[a_{\lambda}(b\bar{\circ}^\hbar_{\mu}c)]_{\hbar}&=&{[a_{\lambda}b]_{\hbar}}\bar{\circ}^\hbar_{\lambda+\mu}c+b\bar{\circ}^\hbar_{\mu}[a_{\lambda}c]_{\hbar},\\
	\label{eq:LJd}[a_{\lambda}[b_{\mu}c]_{\hbar}]_{\hbar}&=&[{[a_{\lambda}b]_{\hbar}}_{\lambda+\mu}c]_{\hbar}+[b_{\mu}[a_{\lambda}c]_{\hbar}]_{\hbar}.
\end{eqnarray}
The $ \hbar $-terms on both sides of \eqref{eq:cLr} give the conformal Leibniz rule: \begin{eqnarray*}
	[a_{\lambda}(b\circ_{\mu}c)]=[a_{\lambda}b]\circ_{\lambda+\mu}c+b\circ _{\mu}[a_{\lambda}c].
\end{eqnarray*}
The $ \hbar^{2} $-terms on both sides of \eqref{eq:LJd} give the conformal Jacobi identity:
\begin{eqnarray*}
	{[a_{\lambda}[b_{\mu}c]]}&=&[[a_{\lambda}b]_{\lambda+\mu}c]+[b_{\mu}[a_{\lambda}c]].
\end{eqnarray*}
By \eqref{eq:semilinear}, we have
$$[\partial a_\lambda b]=-\lambda[a_\lambda b],\quad [a_\lambda \partial b]=(\lambda+\partial)[a_\lambda b],\quad\forall~a,b\in A.$$
Thus $ (A,\circ_\lambda,[\cdot_\lambda\cdot]) $ is a Poisson conformal algebra.
\end{proof}

\begin{defi}
	Let $ (A, \circ_\lambda) $ be a commutative associative conformal algebra. A {\bf conformal $ n $-deformation} of $ A $ is a sequence of bilinear map $ \mu_{i}: A\times A\rightarrow A[\lambda] $, for $ ~0\leqslant i\leqslant n $ with $ \mu_{0} $ being the commutative associative conformal algebra $\lambda$-operation $ \circ_\lambda $ on $ A $ and $\mu_i$ satisfying \eqref{eq:conformal seq} for $ 1\leq i\leq n $ , such that the $ \mathbb{C}[[\hbar]]/\hbar^{n+1} \mathbb{C}[[\hbar]] $-bilinear $\lambda$-operation $\circ^\hbar_\lambda $ on $ A[[\hbar]]/\hbar^{n+1} A[[\hbar]]$ determined by
	\begin{eqnarray*}
		a\circ^\hbar_{\lambda}b=\sum_{k=0}^{n}\hbar^{k}\{a_{\lambda}b\}_{\mu_{k}}=a\circ_{\lambda}b+\sum_{k=1}^{n}\hbar^{k}\{a_{\lambda}b\}_{\mu_{k}},\quad \forall~a, b\in A
	\end{eqnarray*}
	is an associative conformal algebra $\lambda$-operation.
\end{defi}

A conformal $ 1 $-deformation of a commutative associative conformal algebra $ (A, \circ_\lambda) $ is called a {\bf conformal infinitesimal deformation} and we denote it by $ (A, \mu_{1}) $.

Note that  $\mu_1\in C^2_{\rm H}(A,A)$. Furthermore, by direct calculations, $ (A, \mu_{1}) $ is a conformal infinitesimal deformation of a commutative associative conformal algebra $ (A, \circ_\lambda) $ if and only if for all $ a, b, c\in A $
\begin{eqnarray}\label{1defor}
	\{a_{\lambda}(b\circ_{\mu}c)\}_{\mu_{1}}+a\circ_{\lambda}\{b_{\mu}c\}_{\mu_{1}}=\{(a\circ_{\lambda}b)_{\lambda+\mu}c\}_{\mu_{1}}+\{a_{\lambda}b\}_{\mu_{1}}\circ_{\lambda+\mu}c.
\end{eqnarray}
Then we have
\begin{eqnarray*}
	\{a_{\lambda}b_{\mu}c\}_{\dM_{\rm H}\mu_1}
	&=&\{a_{\lambda}(b\circ_{\mu}c)\}_{\mu_{1}}+a\circ_{\lambda}\{b_{\mu}c\}_{\mu_{1}}-\{(a\circ_{\lambda}b)_{\lambda+\mu}c\}_{\mu_{1}}-\{a_{\lambda}b\}_{\mu_{1}}\circ_{\lambda+\mu}c=0,
\end{eqnarray*}
which implies  that $ \mu_{1} $ is a $ 2 $-cocycle for the associative conformal algebra $ (A, \circ_\lambda) $.

Two conformal infinitesimal deformations $ A_{\hbar}=(A, \mu_{1}) $ and $ A^{'}_{\hbar}=(A, \mu^{'}_{1}) $ of a commutative associative conformal algebra $ (A, \circ_\lambda) $ are said to be {\bf equivalent} if there exists a family of conformal algebra homomorphisms ${\rm Id}+\hbar\varphi: A_{\hbar}\rightarrow A^{'}_{\hbar} $ modulo $ \hbar^{2} $.
A conformal infinitesimal deformation is said to be {\bf trivial} if there exists a family of associative conformal  algebra homomorphisms $ {\rm Id}+\hbar\varphi: A_{\hbar}\rightarrow (A, \circ_\lambda) $ modulo $ \hbar^{2} $.

Note that  $\varphi\in C^1_{\rm H}(A,A)$. Furthermore, by direct calculations, $ A_{h} $ and $ A^{'}_{h} $ are equivalent conformal infinitesimal deformations if and only if
\begin{eqnarray}\label{cid}
	\{a_{\lambda}b\}_{\mu_{1}}-\{a_{\lambda}b\}_{\mu^{'}_{1}}=a\circ_{\lambda}\varphi(b)+\varphi(a)\circ_{\lambda}b-\varphi(a\circ_{\lambda}b).
\end{eqnarray}
which implies that $ \mu_{1}-\mu^{'}_{1}=\dM_{\rm H}\varphi $, that is, $ \mu_{1}-\mu^{'}_{1} $ is exact.

Thus we have the following theorem.
\begin{thm}
	There is a one-to-one correspondence between the space of equivalence classes of conformal infinitesimal deformations of $ A $ and the second cohomology group $ H_{\rm H}^{2}(A, A) $.
\end{thm}

It is routine to check that
\begin{pro}
	Let $ (A, \circ_{\lambda}) $ be a commutative associative conformal algebra such that  $ H_{\rm H}^{2}(A, A)=0 $. Then all conformal infinitesimal deformations of $ A $ are trivial.
\end{pro}

\begin{defi}
	Let $ \{\mu_{1}, \dots, \mu_{n}\} $ be a conformal $ n $-deformation of a commutative associative conformal algebra $ (A, \circ_{\lambda}) $. A conformal $ (n+1) $-deformation of a commutative associative conformal algebra $ (A, \circ_{\lambda}) $ given by $ \{\mu_{1}, \dots, \mu_{n}, \mu_{n+1}\} $ is called an {\bf extension} of the conformal $ n $-deformation given by $ \{\mu_{1}, \dots, \mu_{n}\} $.
\end{defi}

\begin{thm}
	For any conformal $ n $-deformation of a commutative associative conformal algebra $ (A, \circ_{\lambda}) $, the map $ \theta_{n}\in C_{\rm H}^{3}(A,A)$ defined by
	\begin{eqnarray}
		\{a_{\lambda}b_{\mu}c\}_{\theta_{n}}=\sum_{r+s=n+1,r, s\geqslant 1}(\{{\{a_{\lambda}b\}_{\mu_{s}}}_{\lambda+\mu}c\}_{\mu_{r}}-\{a_{\lambda}\{b_{\mu}c\}_{\mu_{s}}\}_{\mu_{r}}),\quad \forall a, b, c\in A,
	\end{eqnarray}
    is a $3$-cocycle, i.e. $ \dM_{\rm H}\theta_{n}=0 $.

    Moreover, the conformal $n$-deformation $ \{\mu_{1}, \dots, \mu_{n}\} $ extends into some conformal $ (n+1) $-deformation if and only if $ ~[\theta_{n}]=0 $ in $ H_{\rm H}^{3}(A, A) $.
\end{thm}
\begin{proof}
	It is straightforward to check that $\theta_{n}$ satisfies \eqref{eq:deformCA} and \eqref{eq:deforLA}. That is, $ \theta_{n}\in C_{\rm H}^{3}(A,A)$.
	
	By \eqref{eq:d}, we have
	\begin{eqnarray}
		\nonumber&&\{{a_{1}}_{\lambda_{1}}{a_{2}}_{\lambda_{2}}{a_{3}}_{\lambda_{3}}a_{4}\}_{\dM_{\rm H}\theta_{n}}
	\\
 \nonumber 	&=&{a_{1}}\circ_{\lambda_{1}}\{{a_{2}}_{\lambda_{2}}{a_{3}}_{\lambda_{3}}a_{4}\}_{\theta_{n}}
		-\{({a_{1}}\circ_{\lambda_{1}}a_{2})_{\lambda_{1}+\lambda_{2}}{a_{3}}_{\lambda_{3}}a_{4}\}_{\theta_{n}}+\{{a_{1}}_{\lambda_{1}}({a_{2}}\circ_{\lambda_{2}}a_{3})_{\lambda_{2}+\lambda_{3}}a_{4}\}_{\theta_{n}}
	\\
	\nonumber	&&-\{{a_{1}}_{\lambda_{1}}{a_{2}}_{\lambda_{2}}({a_{3}}\circ_{\lambda_{3}}a_{4})\}_{\theta_{n}}+{\{{a_{1}}_{\lambda_{1}}{a_{2}}_{\lambda_{2}}a_{3}\}_{\theta_{n}}}\circ_{\lambda^{\flat}_{4}}a_{4}\\
 \nonumber &=&\sum_{r+s=n+1,r, s\geqslant 1}{a_{1}}\circ_{\lambda_{1}}\{{\{{a_{2}}_{\lambda_{2}}a_{3}\}_{\mu_{s}}}_{\lambda_{2}+\lambda_{3}}a_{4}\}_{\mu_{r}}-\sum_{r+s=n+1,r, s\geqslant 1}{a_{1}}\circ_{\lambda_{1}}\{{a_{2}}_{\lambda_{2}}\{{a_{3}}_{\lambda_{3}}a_{4}\}_{\mu_{s}}\}_{\mu_{r}}\\
 \nonumber &&-\sum_{r+s=n+1,r, s\geqslant 1}\{{\{({a_{1}}\circ_{\lambda_{1}}a_{2})_{\lambda_{1}+\lambda_{2}}a_{3}\}_{\mu_{s}}}_{\lambda_{1}+\lambda_{2}+\lambda_{3}}a_{4}\}_{\mu_{r}}+\sum_{r+s=n+1,r, s\geqslant 1}\{({a_{1}}\circ_{\lambda_{1}}a_{2})_{\lambda_{1}+\lambda_{2}}\{{a_{3}}_{\lambda_{3}}a_{4}\}_{\mu_{s}}\}_{\mu_{r}}\\
 \label{eq:cocycle} &&+\sum_{r+s=n+1,r, s\geqslant 1}\{{\{{a_{1}}_{\lambda_{1}}({a_{2}}\circ_{\lambda_{2}}a_{3})\}_{\mu_{s}}}_{\lambda_{1}+\lambda_{2}+\lambda_{3}}a_{4}\}_{\mu_{r}}-\sum_{r+s=n+1,r, s\geqslant 1}\{{a_{1}}_{\lambda_{1}}\{({a_{2}}\circ_{\lambda_{2}}a_{3})_{\lambda_{2}+\lambda_{3}}a_{4}\}_{\mu_{s}}\}_{\mu_{r}}\\
 \nonumber &&-\sum_{r+s=n+1,r, s\geqslant 1}\{{\{{a_{1}}_{\lambda_{1}}a_{2}\}_{\mu_{s}}}_{\lambda_{1}+\lambda_{2}}({a_{3}}\circ_{\lambda_{3}}a_{4})\}_{\mu_{r}}
    +\sum_{r+s=n+1,r, s\geqslant 1}\{{a_{1}}_{\lambda_{1}}\{{a_{2}}_{\lambda_{2}}({a_{3}}\circ_{\lambda_{3}}a_{4})\}_{\mu_{s}}\}_{\mu_{r}}\\
 \nonumber   &&+\sum_{r+s=n+1,r, s\geqslant 1}{\{{\{{a_{1}}_{\lambda_{1}}a_{2}\}_{\mu_{s}}}_{\lambda_{1}+\lambda_{2}}a_{3}\}_{\mu_{r}}}\circ_{\lambda^{\flat}_{4}}a_{4}-\sum_{r+s=n+1,r, s\geqslant 1}{\{{a_{1}}_{\lambda_{1}}\{{a_{2}}_{\lambda_{2}}a_{3}\}_{\mu_{s}}\}_{\mu_{r}}}\circ_{\lambda^{\flat}_{4}}a_{4},
	\end{eqnarray}
	where $ \lambda^{\flat}_{4}=\lambda_1+\lambda_2+\lambda_3 $.
	
By the associativity of associative conformal algebra, the third term is the negative of the fifth term and the sixth term is the negative of the eighth term in \eqref{eq:cocycle}.

	By \eqref{fcd}, we have
	\begin{eqnarray}\label{fcde}
		&&\{a_{\lambda}(b\circ_{\mu}c)\}_{\mu_{n}}
		+a\circ_{\lambda}\{b_{\mu}c\}_{\mu_{n}}
		+\sum_{r+s=n+1,r, s\geqslant 1}\{a_{\lambda}\{b_{\mu}c\}_{\mu_{s}}\}_{\mu_{r}}\\
		\nonumber&=&\{(a\circ_{\lambda}b)_{\lambda+\mu}c\}_{\mu_{n}}
		+{\{a_{\lambda}b\}_{\mu_{n}}}\circ_{\lambda+\mu}c
		+\sum_{r+s=n+1,r, s\geqslant 1}\{{\{a_{\lambda}b\}_{\mu_{s}}}_{\lambda+\mu}c\}_{\mu_{r}},
	\end{eqnarray}
which implies that the sum of the rest six terms in \eqref{eq:cocycle} is equal to zero. Thus, $ \theta_{n} $ is closed.
	
	Assume that the conformal $ (n+1) $-deformation of the commutative associative conformal algebra $ (A, \circ_{\lambda}) $ given by $ \{\mu_{1}, \dots, \mu_{n}, \mu_{n+1}\} $ is an extension of the conformal $ n $-deformation given by $ \{\mu_{1}, \dots, \mu_{n}\} $, then we have
	\begin{eqnarray*}
		&&\{a_{\lambda}(b\circ_{\mu}c)\}_{\mu_{n+1}}
		-\{(a\circ_{\lambda}b)_{\lambda+\mu}c\}_{\mu_{n+1}}
		+a\circ_{\lambda}\{b_{\mu}c\}_{\mu_{n+1}}
		-{\{a_{\lambda}b\}_{\mu_{n+1}}}\circ_{\lambda+\mu}c\\
		&=&\sum_{r+s=n+1,r, s\geqslant 1}\{{\{a_{\lambda}b\}_{\mu_{s}}}_{\lambda+\mu}c\}_{\mu_{r}}
		-\sum_{r+s=n+1,r, s\geqslant 1}\{a_{\lambda}\{b_{\mu}c\}_{\mu_{s}}\}_{\mu_{r}}.
	\end{eqnarray*}
	It is obvious that the right-hand side of the above equality is just $ \{a_{\lambda}b_{\mu}c\}_{\theta_{n}} $. We can rewrite the above equality as
	\begin{eqnarray*}
		\{a_{\lambda}b_{\mu}c\}_{\dM_{\rm H}\mu_{n+1}}=\{a_{\lambda}b_{\mu}c\}_{\theta_{n}}.
	\end{eqnarray*}
	We conclude that, if a  conformal $ n $-deformation of commutative associative conformal algebra $ (A,\circ_{\lambda}) $ extends to a conformal $ (n+1) $-deformation, then $ \theta_{n} $ is a coboundary.
	
	Conversely, if $ \theta_{n} $ is a coboundary, then there exists an element $ \psi\in C_{\rm H}^{2}(A, A) $ such that
	\begin{eqnarray*}
		\{a_{\lambda}b_{\mu}c\}_{\dM_{\rm H}\psi}=\{a_{\lambda}b_{\mu}c\}_{\theta_{n}}.
	\end{eqnarray*}
	It is not hard to check that $ \{\mu_{1}, \dots, \mu_{n}, \mu_{n+1}\} $ with $ \mu_{n+1}=\psi $ generates a conformal $(n+1)$-deformation of $ (A, \circ_{\lambda}) $ and thus this conformal $ (n+1) $-deformation is an extension of the conformal $ n $-deformation given by $ \{\mu_{1}, \dots, \mu_{n}\} $.
\end{proof}

\section{Cohomology and deformations of noncommutative Poisson conformal algebras}\label{sec:cohomology}
In this section, we study the Flato-Gerstenhaber-Voronov cohomology theory for noncommutative Poisson conformal algebra associated to a Poisson conformal  module and then use this cohomology to study linear deformations of noncommutative Poisson conformal algebras.

\subsection{Cohomology of noncommutative Poisson conformal algebras}
First, we recall the notion of modules over Lie conformal algebras.
\begin{defi}
	A {\bf module} $V$ over a Lie conformal algebra $A$ is a $\Comp [\partial]$-module endowed with a $\Comp$-bilinear map $A\times V\rightarrow V[\lambda]$, $(a,v)\rightarrow a_{\lambda}v$, satisfying the following axioms:
	\begin{eqnarray}
		\label{eq:Lie module1}  (\partial a)_{\lambda}v&=&-\lambda a_{\lambda}v,\\
		\label{eq:Lie module2}  {a_{\lambda}(\partial v)}&=&(\partial+\lambda)a_{\lambda}v,\\
		\label{eq:Lie module3}  {[a_{\lambda}b]_{\lambda+\mu}v}&=&a_{\lambda}(b_{\mu}v)-b_{\mu}(a_{\lambda}v),\quad\forall~ a,b\in A,v\in V.
	\end{eqnarray}
\end{defi}

The following conclusion is direct.
\begin{pro}
	Let $A$ be a  Lie conformal algebra and let $V$ be a $\Comp [\partial]$-module endowed with a $\Comp$-bilinear map $A\times V\rightarrow V[\lambda]$, $(a,v)\rightarrow a_{\lambda}v$. Then $V$ is a module over the Lie conformal algebra $A$ if and only if $A\oplus V$
	endowed with a $\Comp[\partial]$-module structure and a $\lambda$-bracket defined by:
	\begin{eqnarray*}
		\partial(a+v)&=&\partial a+\partial v, \\
		{[(a+u)_{\lambda}(b+v)]}&=&[a_{\lambda}b]+a_{\lambda}v-b_{-\lambda-\partial}u ,\quad \forall~a,b\in A,u,v\in V
	\end{eqnarray*}
	is a Lie conformal algebra. This Lie conformal algebra is called the {\bf semi-direct product} of $A$ and $V$.
\end{pro}

Note that the structure of a finite module $V$ over a Lie conformal algebra $A$ is the same as a homomorphism of Lie conformal algebras $\rho:A\rightarrow{\rm gc}(V)$. It is obvious that $(\Comp;\rho=0)$ is a module over the Lie conformal algebra $A$, which we call the {\bf trivial module}.

Let $(V;\rho)$ be a module over a  Lie conformal algebra $A$. Define $\rho^{*}:A\rightarrow {\rm gc}(V^{*c})$  by
\begin{equation*}
	(\rho^{*}(a)_{\lambda}\varphi)_{\mu}u=-\varphi_{\mu-\lambda}(\rho(a)_{\lambda}u),
\end{equation*}
for all $a\in A$, $\varphi\in V^{*c}$, $u\in V$. Then $(V^*;\rho^*)$ is a module over $A$.

Assume that $A$ is a finite Lie conformal algebra. Define $\ad:A\lon{\rm gc}(A)$ by $\ad(a)_\lambda b=[a_{\lambda}b]$ for all $a,b\in A$. Then  $(A;\ad)$ is a module over $A$, which we call the {\bf adjoint module}. Furthermore,  $(A^*;\ad^*)$ is also a module over $A$, which we call the {\bf coadjoint module}.

\begin{defi}
	Let $(P,\circ_\lambda,[\cdot_\lambda \cdot])$ be a noncommutative Poisson conformal algebra. A {\bf module} $ V $ over $P$ is a $ \mathbb{C}[\partial] $-module  endowed with  $ \mathbb{C} $-bilinear maps $ A\times V \rightarrow V[\lambda]: ~(a,v) \mapsto a\circ_{\lambda}v $, $ V\times A \rightarrow V[\lambda]:~(v,a) \mapsto v\circ_{\lambda} a $, and $A\times V\rightarrow V[\lambda]:~(a,v)\rightarrow a_{\lambda}v$ such that  $V$ is both the module over the associative conformal algebra $(P,\circ_\lambda)$ and the module over the Lie conformal algebra $(P,[\cdot_\lambda \cdot])$ satisfying
	\begin{eqnarray}
\label{eq:Poisson module 1}		[a_\lambda b]\circ_{\lambda+\mu}v&=&a_\lambda(b\circ_\mu v)-b\circ_\mu(a_\lambda v),\\
	\label{eq:Poisson module 2}		v\circ_{\mu}[a_\lambda b]&=&a_\lambda(v\circ_{\mu}b)-(a_\lambda v)\circ_{\lambda+\mu}b,\\
	\label{eq:Poisson module 3}		(a\circ_\lambda b)_{-\mu-\partial}v&=&a\circ_\lambda(b_{-\mu-\partial} v)+(a_{-\mu-\partial}v)\circ_{\lambda+\mu}b,
	\end{eqnarray}
	where $a,b\in P$ and $v\in V$.
\end{defi}
For convenience, we also set $l(a)_\lambda v=a\circ_\lambda v$, $r(a)_\lambda v=v\circ_{-\lambda-\partial}a$ and $\rho_\lambda(a)v=a_\lambda v$. We denote this module by $(V;l,r,\rho)$.

\begin{ex}
	Let $(P,\circ_\lambda,[\cdot_\lambda \cdot])$ be a noncommutative Poisson conformal algebra. Then $(V;L,R,\ad)$ is a module over the noncommutative Poisson conformal algebra $P$. In general, $(V^{*c};R^*,L^*,\ad^*)$ is not a module over the $P$. However, if $(P,\circ_\lambda,[\cdot_\lambda \cdot])$ is a Poisson conformal algebra, then $(V^{*c};R^*,L^*,\ad^*)$ is also  a module over the $P$.
\end{ex}

\begin{pro}
	Let $P$ be a noncommutative Poisson conformal algebra and let $V$ be a $\Comp [\partial]$-module endowed with $ \mathbb{C} $-bilinear maps $ P\times V \rightarrow V[\lambda]: ~(a,v) \mapsto a\circ_{\lambda}v $, $ V\times P \rightarrow V[\lambda]:~(v,a) \mapsto v\circ_{\lambda} a $, and $P\times V\rightarrow V[\lambda]:~(a,v)\rightarrow a_{\lambda}v$. Then $V$ is a module over the Poisson conformal algebra $P$ if and only if $P\oplus V$
	endowed with a $\Comp[\partial]$-module structure, a  $\lambda$-operation and a $\lambda$-bracket defined by:
	\begin{eqnarray*}
		\partial(a+v)&=&\partial a+\partial v, \\
		{(a+u)\circ_{\lambda}(b+v)}&=&a\circ_{\lambda}b+a\circ_{\lambda}v+u\circ_{\lambda}b, \\
		{[(a+u)_{\lambda}(b+v)]}&=&[a_{\lambda}b]+a_{\lambda}v-b_{-\lambda-\partial}u ,\quad \forall~a,b\in P,u,v\in V
	\end{eqnarray*}
	is a noncommutative Poisson conformal algebra. This noncommutative Poisson conformal algebra is called the {\bf semi-direct product} of $P$ and $V$.
\end{pro}
\begin{proof}
	It follows by a direct calculation.
\end{proof}

The cohomology complex for a Lie conformal algebra $A$ with a module $V$ in the language of $\lambda$-brackets is given as follows (\cite{DK09}). We let $C_{\rm CE}^{0}(A,V)=V/\partial V$ and for $k\geq1$, we denote by $C_{\rm CE}^{k}(A,V)$ the space of all $k$-$\lambda$-brackets on $A$ with coefficients in $V$ satisfying for any $\varphi\in C_{\rm CE}^{k}(A,V)$,
\begin{eqnarray}
	\{a_{1\lambda_{1}}\cdots a_{k-1 \lambda_{k-1}}a_{k}\}_{\varphi}&=&{\rm sign}(\sigma)\{a_{\sigma(1)_{\lambda_{\sigma(1)}}}\cdots a_{\sigma(k-1)_{\lambda_{\sigma(k-1)}}}a_{\sigma(k)}\}_{\varphi}|_{\lambda_{k}\mapsto\lambda_{k}^{\dag}},
\end{eqnarray}
where the notation $\lambda_{k}\mapsto\lambda_{k}^{\dag}$ means that $\lambda_{k}$ is replaced by $\lambda_{k}^{\dag}=-\sum_{j=1}^{k-1}\lambda_{j}-\partial$, if it occurs, and $\partial$ is moved to the left.

Define $C_{\rm CE}^{\bullet}(A,V)=\oplus _{k\in\Nat}C_{\rm CE}^{k}(A,V)$. The corresponding  conformal Chevalley-Eilenberg coboundary operator $\dM_{\rm CE}:C_{\rm CE}^k(A,V)\longrightarrow C_{\rm CE}^{k+1}(A,V)$ is given by
\begin{eqnarray}
	\nonumber    \{a_{1\lambda_{1}}\cdots a_{k\lambda_{k}}a_{k+1}\}_{\dM_{\rm CE} \varphi}=&&\sum_{i=1}^{k}(-1)^{i+1}a_{i\lambda_{i}}\{a_{1\lambda_{1}}\cdots\hat{a}_{i\lambda_i}\cdots a_{k\lambda_{k}}a_{k+1}\}_{\varphi}\\ \label{eq:Lie cohomology}&&+\sum_{i,j=1,i<j}^{k}(-1)^{k+i+j+1}\{a_{1\lambda_{1}}\cdots\hat{a}_{i\lambda_i}\cdots\hat{a}_{j\lambda_j}\cdots a_{k\lambda_{k}}a_{k+1{\lambda_{k+1}^{+}}}[a_{i\lambda_{i}}a_{j}]\}_{\varphi}\\
	\nonumber      &&+(-1)^{k}a_{k+1{\lambda_{k+1}^{\dag}}}\{a_{1\lambda_{1}}\cdots a_{k-1{\lambda_{k-1}}}a_{k}\}_{\varphi}\\
	\nonumber      &&+\sum_{i=1}^{k}(-1)^{i}\{a_{1\lambda_{1}}\cdots\hat{a}_{i\lambda_i}\cdots a_{k\lambda_{k}}[a_{i\lambda_{i}}a_{k+1}]\}_{\varphi},
\end{eqnarray}
where $\varphi\in C_{\rm CE}^k(A,V) $, $\lambda_{k+1}^{\dag}=-\sum_{j=1}^{k}\lambda_{j}-\partial$ and $a_1,\cdots,a_{k+1}\in A$. We denote by $H_{\rm CE}^{\bullet}(A,V)=\oplus _{k\in\Nat}H_{\rm CE}^{k}(A,V)$ the cohomology of the Lie conformal algebra $A$ with coefficients in $V$.

\emptycomment{For $x\in P$ and $\gamma\in C_{\rm H}^{n}(P,V)$, define
	\begin{eqnarray}
	&&	\{a_{1\lambda_{1}}\cdots a_{n \lambda_{n}}a_{n+1}\}_{ \dM_H (x_\lambda \gamma)}\\
	\nonumber			&=&{a_{1}}\circ_{\lambda_{1}}\{{a_{2}}_{\lambda_{2}}\dots {a_{n}}_{\lambda_{n}}{ a_{n+1}}\}_{x_\lambda\gamma}\\
	\nonumber			&&+\sum_{i=1}^{n}(-1)^{i} \{{a_{1}}_{\lambda_{1}}\dots {a_{i-1}}_{\lambda_{i-1}}({a_{i}}\circ_{\lambda_{i}} a_{i+1})_{\lambda_{i}+\lambda_{i+1}}{a_{i+2}}_{\lambda_{i+2}}\dots {a_{n}}_{\lambda_{n}}{a_{n+1}}\}_{x_\lambda\gamma}\\
	%\nonumber &&+(-1)^n\{{a_{1}}_{\lambda_{1}}\dots a_{n-1\lambda_{n-1}}({a_{n}}\circ_{\lambda_{n}}{a_{n+1}})\}_{x_\lambda\gamma}\\
	\nonumber&&+(-1)^{n+1}{\{{a_{1}}_{\lambda_{1}}\dots {a_{n-1}}_{\lambda_{n-1}}{ a_{n}}\}_{x_\lambda\gamma}}\circ_{\tilde{\lambda}^\flat_{n+1}} a_{n+1},
\end{eqnarray}
where ${\tilde{\lambda}^\flat_{n+1}}=\lambda_1+\cdots+\lambda_n+\lambda$.

\begin{pro}
	Let $P$ be a noncommutative Poisson conformal algebra and $V$ a module over $P$. Then  $\dM^2_H (x_\lambda \gamma)=0$.
	Furthermore, we have
	$$x_\lambda \dM_H \gamma=\dM_H (x_\lambda \gamma). $$
\end{pro}
\begin{proof}
On the one hand, by a direct calculation, we have
\begin{eqnarray*}
&&	\{a_{1\lambda_{1}}\cdots a_{n \lambda_{n}}a_{n+1}\}_{x_\lambda \dM_H \gamma}\\
	&=&x_\lambda \{a_{1\lambda_{1}}\cdots a_{n \lambda_{n}}a_{n+1}\}_{\dM_H\gamma}-\sum_{i=1}^{n+1}\{a_{1\lambda_{1}}\cdots([x_\lambda a_i])_{\lambda+\lambda_i}\cdots a_{n \lambda_{n}}a_{n+1}\}_{\dM_H\gamma}
\end{eqnarray*}
\end{proof}}

Let $P$ be a noncommutative Poisson conformal algebra and $V$ a module over $P$. We denote the $(m+n)$-$\lambda$-bracket on $P$ with coefficients in $V$ by
$\gamma:A^{\otimes m}\otimes A^{\otimes n}\rightarrow \Comp[\lambda_{1},\cdots,\lambda_{m},\mu_1,\cdots \mu_{n-1}]\otimes V$ with
$$x_{1}\otimes\cdots\otimes x_{m}\otimes  a_{1 }\otimes \cdots \cdots\otimes  a_{n}\mapsto \{x_{1\lambda_{1}}\cdots x_{m \lambda_{m}}a_{1\mu_1}\cdots a_{n-1\mu_{n-1}}a_n\}_{\gamma}.$$

An $(m+0)$-$\lambda$-bracket $\gamma$ on $P$ with coefficients in $V$  is called an $(m,0)$-cochain if $\gamma$ satisfies
\begin{eqnarray*}
	\{x_{1\lambda_{1}}\cdots x_{m-1 \lambda_{m-1}}x_{m}\}_{\gamma}&=&{\rm sign}(\sigma)\{x_{\sigma(1){\lambda_{\sigma(1)}}}\cdots x_{\sigma(m-1){\lambda_{\sigma(m-1)}}}x_{\sigma(m)}\}_{\gamma}|_{\lambda_{m}\mapsto\lambda_{m}^{\dag}},
\end{eqnarray*}
where the notation $\lambda_{m}\mapsto\lambda_{m}^{\dag}$ means that $\lambda_{m}$ is replaced by $\lambda_{m}^{\dag}=-\sum_{j=1}^{m-1}\lambda_{j}-\partial$.

For $n\geq 1$, an $(m+n)$-$\lambda$-bracket $\gamma$ on $P$ with coefficients in $V$  is called an $(m,n)$-cochain if  $\gamma$ satisfies
\begin{eqnarray*}
	\{x_{1\lambda_{1}}\cdots x_{m \lambda_{m}}a_{1\mu_1}\cdots a_{n-1\mu_{n-1}}a_n\}_{\gamma}= {\rm sign}(\sigma)\{x_{\sigma(1){\lambda_{\sigma(1)}}}\cdots x_{\sigma(m){\lambda_{\sigma(m)}}}a_{1\mu_1}\cdots a_{n-1\mu_{n-1}}a_n\}_{\gamma},
\end{eqnarray*}
where $\sigma$ is a permutation of $\{1,2,\cdots,m\}$.

Denote by $C^{m,n}(P,V)$ the space of $(m,n)$-cochains.
Denote by  $C_{\FGV}^k(P,V)=\sum\limits_{m+n=k\atop m\neq 1}C^{m,n}(P,V)$, which is the space of $k$-cochains.

The {\bf conformal Flato-Gerstenhaber-Voronov coboundary operator} $\dM_{\FGV}:C_{\FGV}^k(P,V)\longrightarrow C_{\FGV}^{k+1}(P,V)$ is defined by
$$\dM_{\FGV}=\sum_{m+n=k}(\dM_{\CE}^{m,n}+(-1)^m\dM_{\rm H}^{m,n} ),$$
where $\dM^{m,n}_{\CE}:C^{m,n}(P,V)\longrightarrow C^{m+1,n}(P,V)$ for $n=0$ is given by \eqref{eq:Lie cohomology} and for $n\geq 1$ is defined by
\begin{eqnarray}\label{eq:FGV-Lie}
	\nonumber   &&	\{x_{1\lambda_{1}}\cdots x_{m+1 \lambda_{m+1}}a_{1\mu_1}\cdots a_{n-1\mu_{n-1}}a_n\}_{\dM^{m,n}_{\rm CE} \gamma}\\
\nonumber	&=&\sum_{i=1}^{m+1}(-1)^{i+1}\Big(x_{i\lambda_{i}}\{x_{1\lambda_{1}}\cdots\hat{x_{i\lambda_{i}}}\cdots x_{m+1 \lambda_{m+1}}a_{1\mu_1}\cdots a_{n-1\mu_{n-1}}a_n\}_{ \gamma}\\
&&-\sum_{j=1}^n\{x_{1\lambda_{1}}\cdots\hat{x_{i\lambda_{i}}}\cdots x_{m+1 \lambda_{m+1}}a_{1\mu_{1}}\cdots([x_{i\lambda_i} a_j])_{\lambda_i+\mu_j}\cdots a_{n-1 \mu_{n-1}}a_{n}\}_{\gamma}
\Big)\\
\nonumber &&+\sum_{i,j=1,i<j}^{m+1}(-1)^{i+j}\{([x_{i\lambda_i}x_j])_{\lambda_i+\lambda_j}x_{1\lambda_{1}}\cdots\hat{x_{i\lambda_i}}\cdots\hat{x_{j\lambda_j}}\cdots x_{m+1 \lambda_{m+1}}a_{1\mu_1}\cdots a_{n-1\mu_{n-1}}a_n\}_{\gamma},
\end{eqnarray}
 $\dM^{m,n}_{\rm H}:C^{m,n}(P,V)\longrightarrow C^{m,n+1}(P,V)$ for $n\geq 2$ is defined by the conformal  Hochschild coboundary operator
\begin{eqnarray*}
	\nonumber   &&	\{x_{1\lambda_{1}}\cdots x_{m \lambda_{m}}a_{1\mu_1}\cdots a_{n\mu_{n}}a_{n+1}\}_{\dM^{m,n}_{\rm H} \gamma}\\
&=&	{a_{1}}\circ_{\mu_{1}}\{x_{1\lambda_{1}}\cdots x_{m \lambda_{m}}a_{2\mu_2}\cdots a_{n\mu_{n}}a_{n+1}\}_{\gamma}\\
	&&+\sum_{i=1}^{n}(-1)^{i} \{{x_{1\lambda_{1}}\cdots x_{m \lambda_{m}}a_{1}}_{\mu_{1}}\dots {a_{i-1}}_{\mu_{i-1}}({a_{i}}\circ_{\mu_{i}} a_{i+1})_{\mu_{i}+\mu_{i+1}}{a_{i+2}}_{\mu_{i+2}}\dots {a_{n}}_{\mu_{n}}{a_{n+1}}\}_{\gamma}\\
	%\nonumber &&+(-1)^n\{{a_{1}}_{\lambda_{1}}\dots a_{n-1\lambda_{n-1}}({a_{n}}\circ_{\lambda_{n}}{a_{n+1}})\}_{\gamma}\\
	\nonumber&&+(-1)^{n+1}{\{x_{1\lambda_{1}}\cdots x_{m \lambda_{m}}{a_{1}}_{\mu_{1}}\dots {a_{n-1}}_{\mu_{n-1}}{ a_{n}}\}_{\gamma}}\circ_{\mu^\flat_{n+1}} a_{n+1},
\end{eqnarray*}
where $\mu^\flat_{n+1}=\sum_{i=1}^m\lambda_i+\sum_{j=1}^n\mu_j$,  and $\dM^{m-1,1}_{\rm H}$ is the composition of the natural {linear map}
$C^{m,0}(P,V)\hookrightarrow C^{m-1,1}(P,V)$
and the conformal Hochschild coboundary operator $\dM_{\rm H}:C^{m-1,1}(P,V)\longrightarrow C^{m-1,2}(P,V)$ given by
\begin{eqnarray*}
	\nonumber   	\{x_{1\lambda_{1}}\cdots x_{m \lambda_{m}}a_{1}\}_{\dM_{\rm H} \gamma}&=&	{x_{m}}\circ_{\lambda_{m}}\{x_{1\lambda_{1}}\cdots x_{m-1 \lambda_{m-1}}a_{1}\}_{\gamma}-\{x_{1\lambda_{1}}\cdots x_{m-1 \lambda_{m-1}} (x_{m}\circ_{\lambda_{m}} a_1)\}_{\gamma}\\
	&&+\{x_{1\lambda_{1}}\cdots x_{m-1\lambda_{m-1}}x_{m }\}_{\gamma}\circ_{\mu^\flat_{1}}a_1,
\end{eqnarray*}
where $\mu^\flat_{1}=\sum_{i=1}^m\lambda_i$.

It is illustrated by the following diagram of bicomplex:
{\begin{equation*}\label{diagram1}
	\begin{array}{ccccccccc}
		&            &\cdots&       &\cdots&          &\cdots& &\\
		&            &\Big\uparrow&       &\Big\uparrow&          &\Big\uparrow& &\\
		0&\longrightarrow& C^{0,3}(P,V)&\stackrel{\dM^{0,3}_{\CE}}\longrightarrow&C^{1,3}(P,V)&\stackrel{\dM^{1,3}_{\CE}}\longrightarrow&C^{2,3}(P,V)&\longrightarrow&\cdots\\
		&            &\dM^{0,2}_{\rm H}\Big\uparrow&       &\dM^{1,2}_{\rm H}\Big\uparrow&          &\dM^{2,2}_{\rm H}\Big\uparrow& &\\
		0&\longrightarrow& C^{0,2}(P,V)&\stackrel{\dM^{0,2}_{\CE}}\longrightarrow&C^{1,2}(P,V)&\stackrel{\dM^{1,2}_{\CE}}\longrightarrow&C^{2,2}(P,V)&\longrightarrow&\cdots\\
		&            &\dM^{0,1}_{\rm H}\Big\uparrow&       &\dM^{1,1}_{\rm H}\Big\uparrow&          &\dM^{2,1}_{\rm H}\Big\uparrow& &\\
		V/\partial V&\stackrel{\dM_{\CE}}\longrightarrow& C^{1,0}(P,V)&\stackrel{\dM^{1,0}_{\CE}}\longrightarrow&C^{2,0}(P,V)&\stackrel{\dM^{2,0}_{\CE}}\longrightarrow&C^{3,0}(P,V)&\longrightarrow&\cdots
\end{array}\end{equation*}}

Let $P$ be a noncommutative Poisson conformal algebra and $V$ a module over $P$.  Define a $\Comp$-linear map $ P\times C_{\rm H}^{n}(P,V) \rightarrow C_{\rm H}^{n}(P,V)[\lambda]: ~(x,\gamma) \mapsto x_{\lambda}\gamma $ by
\begin{align}\label{eq:rep-FGV}
	\{a_{1\lambda_{1}}\cdots a_{n-1 \lambda_{n-1}}a_{n}\}_{x_\lambda\gamma}=x_\lambda \{a_{1\lambda_{1}}\cdots a_{n-1 \lambda_{n-1}}a_{n}\}_{\gamma}-\sum_{i=1}^n\{a_{1\lambda_{1}}\cdots([x_\lambda a_i])_{\lambda+\lambda_i}\cdots a_{n-1 \lambda_{n-1}}a_{n}\}_{\gamma},
\end{align}
where $a_1,\cdots,a_n\in P$.  For $0\leq i\leq n-1$, we have
\begin{align*}
	&	\{a_{1\lambda_{1}}\cdots(\partial a_{i })_{\lambda_i}\cdots a_{n-1 \lambda_{n-1}}a_{n}\}_{x_\lambda\gamma}\\
	&=	x_\lambda 	\{a_{1\lambda_{1}}\cdots(\partial a_{i })_{\lambda_i}\cdots a_{n-1 \lambda_{n-1}}a_{n}\}_{\gamma}-\sum_{j=1}^{i-1}\{a_{1\lambda_{1}}   \cdots([x_\lambda a_j])_{\lambda+\lambda_i}\cdots(\partial a_{i })_{\lambda_i}\cdots a_{n-1 \lambda_{n-1}}a_{n}\}_{\gamma}\\
	&-\{a_{1\lambda_{1}}\cdots([x_\lambda (\partial a_i)])_{\lambda+\lambda_i}\cdots a_{n-1 \lambda_{n-1}}a_{n}\}_{\gamma}-\sum_{j=i+1}^{n}\{a_{1\lambda_{1}}\cdots(\partial a_{i })_{\lambda_i}   \cdots([x_\lambda a_j])_{\lambda+\lambda_i}\cdots a_{n-1 \lambda_{n-1}}a_{n}\}_{\gamma}\\
	&=-\lambda_i x_\lambda 	\{a_{1\lambda_{1}}\cdots a_{i \lambda_i}\cdots a_{n-1 \lambda_{n-1}}a_{n}\}_{\gamma}+\lambda_i\sum_{j=1}^{i-1}\{a_{1\lambda_{1}}  \cdots([x_\lambda a_j])_{\lambda+\lambda_i}\cdots a_{i \lambda_i} \cdots a_{n-1 \lambda_{n-1}}a_{n}\}_{\gamma}\\
	&-\lambda_i\{a_{1\lambda_{1}}\cdots([x_\lambda a_i])_{\lambda+\lambda_i}\cdots a_{n-1 \lambda_{n-1}}a_{n}\}_{\gamma}+\lambda_{i}\sum_{j=i+1}^{n}\{a_{1\lambda_{1}}\cdots a_{i\lambda_i}   \cdots([x_\lambda a_j])_{\lambda+\lambda_i}\cdots a_{n-1 \lambda_{n-1}}a_{n}\}_{\gamma}\\
	&=-\lambda_i \{a_{1\lambda_{1}}\cdots a_{n-1 \lambda_{n-1}}a_{n}\}_{x_\lambda\gamma},
\end{align*}	
which implies that
\begin{equation}\label{eq:property 1}
	\{a_{1\lambda_{1}}\cdots(\partial a_{i })_{\lambda_i}\cdots a_{n-1 \lambda_{n-1}}a_{n}\}_{x_\lambda\gamma}=	-\lambda_i \{a_{1\lambda_{1}}\cdots a_{n-1 \lambda_{n-1}}a_{n}\}_{x_\lambda\gamma}.
\end{equation}	

Similarly, we also have
\begin{equation}
	\{a_{1\lambda_{1}}\cdots a_{n-1 \lambda_{n-1}}\partial a_{n}\}_{x_\lambda\gamma}=(\lambda_1+\cdots+\lambda_{n-1}+\lambda+\partial)\{a_{1\lambda_{1}}\cdots a_{n-1 \lambda_{n-1}}a_{n}\}_{x_\lambda\gamma}.
\end{equation}
Note that the space  $C_{\rm H}^{n}(P,V)$ endows with a natural $\Comp[\tilde{\partial}]$-module structure defined by
\begin{equation}
	\{a_{1\lambda_{1}}\cdots a_{n-1 \lambda_{n-1}}a_{n}\}_{\tilde{\partial}\gamma}=(\lambda_1+\cdots+\lambda_{n-1}+\partial)\{a_{1\lambda_{1}}\cdots a_{n-1 \lambda_{n-1}}a_{n}\}_{\gamma}.
\end{equation}	

Furthermore,  we have
\begin{pro}\label{pro:rep-FGV}
	For $n\geq 1$, the  $\Comp$-linear map $ P\times C_{\rm H}^{n}(P,V) \rightarrow C_{\rm H}^{n}(P,V)[\lambda]: ~(x,\gamma) \mapsto x_{\lambda}\gamma $ given by \eqref{eq:rep-FGV} defines a module  $C_{\rm H}^{n}(P,V)$  over  the Lie conformal algebra $(P,[\cdot_\lambda\cdot])$.
\end{pro}
\begin{proof}
	By \eqref{eq:conformal seq} and \eqref{eq:Lie module1}, we have
	$$\{a_{1\lambda_{1}}\cdots a_{n-1 \lambda_{n-1}}a_{n}\}_{(\partial x)_\lambda\gamma}=-\lambda\{a_{1\lambda_{1}}\cdots a_{n-1 \lambda_{n-1}}a_{n}\}_{x_\lambda\gamma}.$$
	By \eqref{eq:conformal seq} and \eqref{eq:Lie module2}, we have
	$$\{a_{1\lambda_{1}}\cdots a_{n-1 \lambda_{n-1}}a_{n}\}_{ x_\lambda(\tilde{\partial}\gamma)}=(\lambda+ \lambda_1+\cdots+\lambda_{n-1}+\partial)\{a_{1\lambda_{1}}\cdots a_{n-1 \lambda_{n-1}}a_{n}\}_{x_\lambda\gamma}.$$
	By \eqref{eq:Lie conformal3} and \eqref{eq:Lie module3}, we have
	\begin{eqnarray*}
		&&	\{a_{1\lambda_{1}}\cdots a_{n-1 \lambda_{n-1}}a_{n}\}_{([ x_\lambda y]_{\lambda+\mu} \gamma-x_\lambda (y_\mu \gamma)+y_\mu(x_\lambda \gamma))}\\
		&=&([ x_\lambda y]_{\lambda+\mu})\{a_{1\lambda_{1}}\cdots a_{n-1 \lambda_{n-1}}a_{n}\}_\gamma-x_\lambda (y_\mu \{a_{1\lambda_{1}}\cdots a_{n-1 \lambda_{n-1}}a_{n}\}_\gamma)+y_\mu(x_\lambda \{a_{1\lambda_{1}}\cdots a_{n-1 \lambda_{n-1}}a_{n}\}_\gamma)\\
		&&-\sum_{i=1}^{n}\{a_{1\lambda_{1}}\cdots ([[x_\lambda y]_{\lambda+\mu} a_i]-[x_\lambda [y_\mu a_i]]+[y_\mu[x_\lambda a_i]])_{\lambda+\mu+\lambda_i}\cdots a_{n-1 \lambda_{n-1}}a_{n}\}_{\gamma}\\
		&=&0.
	\end{eqnarray*}
Thus  $C_{\rm H}^{n}(P,V)$  is a module structure over  the Lie conformal algebra $(P,[\cdot_\lambda\cdot])$.
\end{proof}

Recall that the basic cohomology complex   for a Lie conformal algebra $(A, [\cdot_\lambda\cdot])$ associated to a module $V$ is given as follows. An $m$-cochain of a Lie conformal algebra $A$ with coefficients in a module $V$ in the basic cohomology complex is a $\Comp$-linear map
$$\gamma:A^{\otimes m}\rightarrow V[\lambda_1,\cdots,\lambda_{m}],\quad x_1\otimes \cdots\otimes x_m\mapsto \gamma_{\lambda_1,\cdots,\lambda_m}(x_1,\cdots,x_m) $$
satisfying the following conditions:
\begin{eqnarray}
\label{eq:basic cochain 1}	\gamma_{\lambda_1,\cdots,\lambda_m}(x_1,\cdots,\partial x_i,\cdots,x_m)&=&-\lambda_i 	\gamma_{\lambda_1,\cdots,\lambda_m}(x_1,\cdots, x_i,\cdots,x_m),\quad  1\leq i\leq m;\\
\label{eq:basic cochain 2}		\gamma_{\lambda_1,\cdots,\lambda_i,\lambda_{i+1},\cdots,\lambda_m}(x_1,\cdots, x_i,x_{i+1},\cdots,x_m)&=&-	\gamma_{\lambda_1,\cdots,\lambda_{i+1},\lambda_{i},\cdots,\lambda_m}(x_1,\cdots, x_{i+1},x_{i},\cdots,x_m).
\end{eqnarray}

Denote by $\tilde{C}_{\CE}^n(A,V)$ the space of $n$-cochains  the basic cohomology complex. The corresponding  coboundary operator
$\tilde{\dM}_{\CE}:\tilde{C}_{\CE}^m(A,V)\longrightarrow \tilde{C}_{\CE}^{m+1}(A,V)$ is given by
\begin{eqnarray*}
&&(	\tilde{\dM}_{\CE}\gamma)_{\lambda_1,\cdots,\lambda_{m+1}}(x_1,\cdots,x_{m+1})\\
&=&\sum_{i=1}^{m+1}(-1)^{i+1}x_{i\lambda_i}\gamma_{\lambda_1,\cdots,\hat{\lambda_i},\cdots,\lambda_{m+1}}(x_1,\cdots, \hat{x_i},\cdots,x_{m+1})\\
&&+\sum_{i,j=1,i<j,}^{m+1}(-1)^{i+j}\gamma_{\lambda_i+\lambda_j,\lambda_1,\cdots,\hat{\lambda_i},\cdots,\hat{\lambda_i},\cdots,\lambda_{m+1}}
([x_{i\lambda_i}x_j],x_{1},\cdots,\hat{x_{i}},\cdots,\hat{x_{j}},\cdots, x_{m+1}).
	\end{eqnarray*}

\begin{lem}\label{lem:FGV1}
	Let $P$ be a noncommutative Poisson conformal algebra and $V$ a module over $P$.
	\begin{itemize}
		\item[{\rm (i)}] For $\gamma\in C^{m,n}(P,V)$, we have $\dM^{m,n}_{\rm CE} \gamma\in C^{m+1,n}(P,V)$ and $\dM^{m+1,n}_{\rm CE}\circ \dM^{m,n}_{\rm CE}=0$.
		\item[{\rm (ii)}]For $\gamma\in C^{m,n}(P,V)$, we have $\dM^{m,n}_{\rm H} \gamma\in C^{m,n+1}(P,V)$ and $\dM^{m,n+1}_{\rm H} \circ \dM^{m,n}_{\rm H} =0$.
	\end{itemize}
\end{lem}	
	\begin{proof}
(i) For any $\gamma\in C^{m,n}(P,V)$ with $n=0$, this case is just the cohomology of Lie conformal algebra. The conclusions follow immediately. For any $\gamma\in C^{m,n}(P,V)$ with $n\geq 1$, define a $\Comp$-linear map
$$\tilde{\gamma}:A^{\otimes m}\rightarrow C_{\rm H}^{n}(P,V),\quad x_1\otimes \cdots \otimes x_m\mapsto \tilde{\gamma}_{\lambda_1,\cdots,\lambda_m}(x_1,\cdots,x_m)  $$
by
$$\{a_{1\mu_1}\cdots a_{n-1\mu_{n-1}}a_n\}_{\tilde{\gamma}_{\lambda_1,\cdots,\lambda_m}(x_1,\cdots,x_m)  }=	\{x_{1\lambda_{1}}\cdots x_{m \lambda_{m}}a_{1\mu_1}\cdots a_{n-1\mu_{n-1}}a_n\}_{\gamma}.$$
It is obvious that this map is an isomorphism.  By a direct calculation, $\tilde{\gamma}$ satisfies \eqref{eq:basic cochain 1} and \eqref{eq:basic cochain 2}. Thus we have  $\tilde{\gamma}\in \tilde{C}_{\CE}^m(A,C_{\rm H}^{n}(P,V))$. Furthermore, by Proposition \ref{pro:rep-FGV}, the map  defined  by \eqref{eq:rep-FGV} gives a module $C_{\rm H}^{n}(P,V)$ over the Lie conformal algebra $(P,[\cdot_\lambda\cdot])$. We obtain a coboundary operator $\tilde{\dM}_{\CE}:\tilde{C}_{\CE}^m(A,C_{\rm H}^{n}(P,V))\longrightarrow \tilde{C}_{\CE}^{m+1}(A,C_{\rm H}^{n}(P,V))$ given by
\begin{eqnarray*}
&&{\{a_{1\mu_1}\cdots a_{n-1\mu_{n-1}}a_n\}}_{(\tilde{\dM}_{\CE}\tilde{\gamma})_{\lambda_1,\cdots,\lambda_{m+1}}(x_1,\cdots,x_{m+1})  }\\
&=&{\sum_{i=1}^{m+1}(-1)^{i+1}{\{a_{1\mu_1}\cdots a_{n-1\mu_{n-1}}a_n\}}_{x_{i\lambda_i}\tilde{\gamma}_{\lambda_1,\cdots,\hat{\lambda_i},\cdots,\lambda_m}(x_1,\cdots, \hat{x_i},\cdots,x_{m+1})}}\\
&&+\sum_{i,j=1,i<j,}^{m+1}(-1)^{i+j}  \{a_{1\mu_1}\cdots a_{n-1\mu_{n-1}}a_n\}_{ \tilde{\gamma}_{\lambda_i+\lambda_j,\lambda_1,\cdots,\hat{\lambda_i},\cdots,\hat{\lambda_i},\cdots,\lambda_{m+1}}
([x_{i\lambda_i}x_j],x_{1},\cdots,\hat{x_{i}},\cdots,\hat{x_{j}},\cdots, x_{m+1})}\\
&=&\sum_{i=1}^{m+1}(-1)^{i+1}\Big(x_{i\lambda_{i}}\{x_{1\lambda_{1}}\cdots\hat{x_{i\lambda_{i}}}\cdots x_{m+1 \lambda_{m+1}}a_{1\mu_1}\cdots a_{n-1\mu_{n-1}}a_n\}_{ \gamma}\\
&&-\sum_{j=1}^n\{x_{1\lambda_{1}}\cdots\hat{x_{i\lambda_{i}}}\cdots x_{m+1 \lambda_{m+1}}a_{1\mu_{1}}\cdots([x_{i\lambda_i} a_j])_{\lambda_i+\mu_j}\cdots a_{n-1 \lambda_{n-1}}a_{n}\}_{\gamma}
\Big)\\
\nonumber &&+\sum_{i,j=1,i<j}^{m+1}(-1)^{i+j}\{([x_{i\lambda_i}x_j])_{\lambda_i+\lambda_j}x_{1\lambda_{1}}\cdots\hat{x_{i\lambda_i}}\cdots\hat{x_{j\lambda_j}}\cdots x_{m+1 \lambda_{m+1}}a_{1\mu_1}\cdots a_{n-1\mu_{n-1}}a_n\}_{\gamma}\\
&=&\{x_{1\lambda_{1}}\cdots x_{m+1 \lambda_{m+1}}a_{1\mu_1}\cdots a_{n-1\mu_{n-1}}a_n\}_{\dM^{m,n}_{\rm CE} \gamma},
\end{eqnarray*}	
which implies that $\dM^{m,n}_{\rm CE} \gamma\in C^{m+1,n}(P,V)$ and $\dM^{m+1,n}_{\rm CE}\circ \dM^{m,n}_{\rm CE}=0$.

(ii) By the properties of $\lambda$-bracket on $P$ with coefficients in $V$, it is straightforward to check that $\dM^{m,n}_{\rm H} \gamma\in C^{m,n+1}(P,V)$. By the cohomology of the associative conformal algebra $(P,\circ_\lambda)$, we have $\dM^{m,n+1}_{\rm H} \circ \dM^{m,n}_{\rm H} =0$.
\end{proof}
	
\begin{lem}\label{lem:FGV2}
	Let $P$ be a noncommutative Poisson conformal algebra and $V$ a module over $P$. Then we have
	\begin{eqnarray}
		\dM^{m,1}_{\rm H}\circ \dM^{m,0}_{\rm CE}&=&\dM^{m-1,2}_{\rm CE}\circ \dM^{m-1,1}_{\rm H},~m\geq1,\\
		\dM^{m+1,n}_{\rm H}\circ \dM^{m,n}_{\rm CE}&=&\dM^{m,n+1}_{\rm CE}\circ \dM^{m,n}_{\rm H},~m\geq0,n\geq 2.
	\end{eqnarray}

\emptycomment{	\begin{eqnarray}
		\dM^{m-1,2}_{\rm H}\circ \dM^{m-1,1}_{\rm CE}&=&-\dM^{m,1}_{\rm CE}\circ \dM^{m,0}_{\rm H},\quad m\geq 1,\\
	\dM^{m+1,n}_{\rm H}\circ \dM^{m,n}_{\rm CE}&=&\dM^{m,n+1}_{\rm CE}\circ \dM^{m,n}_{\rm H}\quad m\geq 0, n\geq 2.
	\end{eqnarray}}
	
\end{lem}
\begin{proof}
	It is given in Appendix.
\end{proof}

\begin{thm}
	With the above notations, we have $\dM_{\FGV}\circ \dM_{\FGV}=0$, i.e. $(C_{\FGV}^\ast(P,V)=\oplus_{k=0}^{\infty}  C_{\FGV}^k(P,V),\dM_{\FGV})$ is a
	cochain complex.
\end{thm}
\begin{proof}
	By Lemma \ref{lem:FGV1}, Lemma \ref{lem:FGV2} and the theory of double bicomplex and total complex, the conclusion follows.
	\end{proof}

\begin{defi}
		Let $P$ be a noncommutative Poisson conformal algebra and $V$ a module over $P$. The cohomology of the cochain complex $\big(\oplus_{k=0}^{+\infty}C_{\FGV}^k(P,V), \dM_{\FGV})$ is called the {\bf noncommutative Poisson conformal complex  } of $P$ with coefficients in $V$. The corresponding $k$-th cohomology group, denoted by $H_{\FGV}^k(P,V)=Z^k_{\FGV}(P,V)/B^k_{\FGV}(P,V)$, is called the {\bf noncommutative Poisson conformal cohomology}.
	\end{defi}
\begin{rmk}
The cohomology of Poisson conformal algebras can also be obtained by replacing the conformal Hochschild cohomology by conformal Harrison cohomology in the bicomplex and then the analogous assertions are still hold in the below.
\end{rmk}
Note that the $0$-th cohomology of noncommutative Poisson conformal algebra $P$ is given by
$$H_{\FGV}^1(P,V)=\{v\in V\mid a_{-\partial }v=0,~\forall~a\in P  \}/\partial V.$$
By a direct calculation, we have
  \begin{eqnarray*}
	Z^1_{\FGV}(P,V)=\{\gamma\in C^{1,0}(P,V)\mid&& \gamma([a_{1\lambda_1} a_2])=a_{1\lambda_1}\gamma(a_2)-a_{2(-\partial-\lambda_1)}\gamma(a_1)\\
	&& \gamma(a_1\circ_{\lambda_1}a_2)=a_1\circ_{\lambda_1}\gamma(a_2)+\gamma(a_1)\circ_{\lambda_1}a_2,\quad\forall~a_1,a_2\in P\},
\end{eqnarray*}
and
\begin{eqnarray*}
B^1_{\FGV}(P,V)=\{\gamma\in C^{1,0}(P,V)\mid \mbox{there exists  $v\in V$ such that } \gamma(a)=a_{-\partial} v,\quad \forall~a\in P\}.
\end{eqnarray*}
Thus $H_{\FGV}^1(P,V)=Z^1_{\FGV}(P,V)/B^1_{\FGV}(P,V)$.

\subsection{Deformations of noncommutative Poisson conformal algebras}
	Let $P$ be a noncommutative Poisson conformal algebra and $\varpi,\omega:P\otimes P\rightarrow P[\lambda]$ two $\Comp$-linear maps. Consider two $t$-parameterized family of bilinear $\lambda$-brackets and $\lambda$-operations
\begin{eqnarray}
{	a\circ_{t\lambda} b}= a\circ_\lambda b+t\{a_\lambda b\}_{\varpi},\quad	{[a_{\lambda}b]_t}=[a_{\lambda}b]+t\{a_{\lambda}b\}_{\omega},\quad\forall~a,b\in P.
\end{eqnarray}
If  $(P,\circ_{t\lambda},[\cdot_{\lambda}\cdot]_t)$ is a noncommutative Poisson conformal algebra for all $t$, we say
that $(\varpi,\omega)$ generates a  {\bf linear
	deformation} of $P$.

By a direct calculation,  we can deduce that $(P,\circ_{t\lambda},[\cdot_{\lambda}\cdot]_t)$ is a noncommutative Poisson conformal algebra for any $t\in \mathbb \Real$ if and only if
{\small\begin{eqnarray}
		\label{eq:ass conformal1}\{\partial a_{\lambda}b\}_{\varpi}&=&-\lambda\{a_{\lambda}b\}_{\varpi},\\
	\label{eq:ass conformal2} \{ a_{\lambda}\partial b\}_{\varpi}&=&(\partial+\lambda)\{a_{\lambda}b\}_{\varpi},\\
	\label{eq:ass conformal3}	\{a_{\lambda}\{b_{\mu}c\}_{\varpi}\}_{\varpi}&=&\{\{a_{\lambda}b\}_{\varpi_{\lambda+\mu}}c\}_{\varpi},\\
	\label{eq:lie conformal1}\{\partial a_{\lambda}b\}_{\omega}&=&-\lambda\{a_{\lambda}b\}_{\omega},\\
	\label{eq:lie conformal2} \{ a_{\lambda}\partial b\}_{\omega}&=&(\partial+\lambda)\{a_{\lambda}b\}_{\omega},\\
	\label{eq:lie conformal3}\{a_{\lambda}b\}_{\omega}&=&-\{b_{-\lambda-\partial}a\}_{\omega},\\
	\label{eq:lie conformal4}	\{a_{\lambda}\{b_{\mu}c\}_{\omega}\}_{\omega}&=&\{b_{\mu}\{a_{\lambda}c\}_{\omega}\}_{\omega}+\{\{a_{\lambda}b\}_{\omega_{\lambda+\mu}}c\}_{\omega},\\
\label{eq:Poisson Leibniz}	\{a_\lambda\{b_\mu c\}_{\varpi}\}_{\omega}&=&\{\{a_\lambda b\}_{\omega_{\lambda+\mu}}c\}_\varpi+\{b_\mu\{a_\lambda c\}_\omega\}_\varpi,\\
		\label{eq:deformation-FGV1}	\{a_{\lambda}(b\circ_{\mu}c)\}_{\varpi}+a\circ_{\lambda}\{b_{\mu}c\}_{\varpi}&=&\{(a\circ_{\lambda}b)_{\lambda+\mu}c\}_{\varpi}+\{a_{\lambda}b\}_{\varpi}\circ_{\lambda+\mu}c,\\
	\label{eq:deformation-FGV2}\{[a_{\lambda}b]_{\lambda+\mu}c\}_{\omega}-\{a_{\lambda}[b_{\mu}c]\}_{\omega}
	+\{b_{\mu}[a_{\lambda}c]\}_{\omega}&=&[a_{\lambda}\{b_{\mu}c\}_{\omega}]-[\{a_{\lambda}b\}_{\omega_{\lambda+\mu}}c]-[b_{\mu}\{a_{\lambda}c\}_{\omega}],\\
	\label{eq:deformation-FGV3}{	[a_\lambda\{b_\mu c\}_\varpi]}-\{[a_\lambda b]_{\lambda+\mu}c\}_\varpi-	\{b_\mu[a_\lambda c]\}_\varpi&=&\{a_\lambda b\}_\omega\circ_{\lambda+\mu} c+b\circ_\mu\{a_\lambda c\}_\omega-\{a_\lambda(b\circ_\mu c)\}_\omega.
\end{eqnarray}}
It is not hard to see that \eqref{eq:ass conformal1}-\eqref{eq:Poisson Leibniz} mean that $(P,\{\cdot_\lambda\cdot\}_\varpi,\{\cdot_\lambda\cdot\}_\omega)$ is a noncommutative Poisson conformal algebra. \eqref{eq:deformation-FGV1} means that $\{a_\lambda b_\mu c\}_{\dM^{0,2}_H\varpi}=0$, \eqref{eq:deformation-FGV2} means that $\{a_\lambda b_\mu c\}_{\dM^{2,0}_\CE \omega}=0$, and \eqref{eq:deformation-FGV3} means that $\{a_\lambda b_\mu c\}_{\dM^{0,2}_\CE \varpi-\dM^{1,1}_H\omega}=0$.

\begin{pro}
	If $(\varpi,\omega)$ generates a  { linear deformation} of the  noncommutative Poisson conformal algebra  $P$, then $\omega+\varpi$ is $2$-cocycle for the noncommutative Poisson conformal algebra  $P$.
\end{pro}
\begin{proof}
	By \eqref{eq:deformation-FGV1}-\eqref{eq:deformation-FGV3}, we deduce that $\dM_{\FGV}(\omega+\varpi)=0$, i.e. $\omega+\varpi$ is $2$-cocycle for the noncommutative Poisson conformal algebra  $P$.
\end{proof}

\begin{defi}
Two linear deformations	 $P_t=(P,\circ_{t\lambda},[\cdot_{\lambda}\cdot]_t)$ and  $P'_t=(P,\circ'_{t\lambda},[\cdot_{\lambda}\cdot]'_t)$  of the noncommutative Poisson conformal algebra $(P,\circ_{\lambda},[\cdot_{\lambda}\cdot])$, which are generated by  $(\varpi,\omega)$ and $(\varpi',\omega')$ respectively,  are said to be {\bf equivalent} if there exists a family of  noncommutative Poisson  conformal algebra
homomorphisms ${\Id}+tN:P_t\longrightarrow P'_t$.

A deformation
is said to be {\bf trivial} if there exist a family of noncommutative Poisson conformal algebra homomorphisms ${\Id}+tN:P_t\longrightarrow (P,\circ_{\lambda},[\cdot_{\lambda}\cdot])$.
\end{defi}

Two linear deformations $P_t$ and $P'_t$ generated by  $(\varpi,\omega)$ and $(\varpi',\omega')$ respectively of a  noncommutative Poisson conformal algebra $(P,\circ_{\lambda},[\cdot_{\lambda}\cdot])$ are equivalent if and only if
\begin{eqnarray}
	N\circ \partial&=&\partial\circ N,\label{eq:equivalent1}\\
	\{a_\lambda b\}_\varpi-\{a_\lambda b\}_{\varpi'}&=&N(a)\circ_\lambda b+a\circ_\lambda N(b)-N(a\circ_\lambda b),\label{eq:equivalent2}\\
	\{a_\lambda b\}_\omega-\{a_\lambda b\}_{\omega'}&=&[N(a)_\lambda b]+[a_\lambda N(b)]-N[a_\lambda b],\label{eq:equivalent3}\\
	N\{a_\lambda b\}_\varpi&=&\{N(a)_\lambda b\}_{\varpi'}+\{a_\lambda N(b)\}_{\varpi'}+N (a)\circ_\lambda N(b),\label{eq:equivalent4}\\
	\{N(a)_\lambda N(b)\}_{\varpi'}&=&0,\label{eq:equivalent5}\\
	N\{a_\lambda b\}_\omega&=&\{N(a)_\lambda b\}_{\omega'}+\{a_\lambda N(b)\}_{\omega'}+[N (a)_\lambda N(b)],\label{eq:equivalent6}\\
	\{N(a)_\lambda N(b)\}_{\omega'}&=&0.\label{eq:equivalent7}
\end{eqnarray}
It is hard to see that \eqref{eq:equivalent1}-\eqref{eq:equivalent3} imply that
$$(\omega+\varpi)-(\omega'+\varpi')=\dM_{\FGV} N.$$

We summarize the above discussion in the following conclusions.
\begin{thm}
	If two linear deformations $P_t$ and $P'_t$ generated by  $(\varpi,\omega)$ and $(\varpi',\omega')$ respectively of a  noncommutative Poisson conformal algebra $P$ are equivalent, then $\omega+\varpi$ and $\omega'+\varpi'$ are in the same cohomology class in $H^2_{\FGV}(P,P)$.
	\end{thm}
	
One can deduce that $P_t$ is a trivial linear deformation if and only if there exists $N:P\rightarrow P$ such that for all $a,b\in P$, we have
\begin{eqnarray}
	N\circ \partial&=&\partial\circ N,\label{Nij1}\\
		\{a_\lambda b\}_\varpi&=&N(a)\circ_\lambda b+a\circ_\lambda N(b)-N(a\circ_\lambda b),\label{Nij2}\\
		N\{a_\lambda b\}_\varpi&=&N (a)\circ_\lambda N(b),\label{Nij3}\\
	\{a_\lambda b\}_\omega&=&[N(a)_\lambda b]+[a_\lambda N(b)]-N[a_\lambda b],\label{Nij4}\\
	N\{a_\lambda b\}_\omega&=&[N (a)_\lambda N(b)].\label{Nij5}
\end{eqnarray}

It follows from \eqref{Nij2} and \eqref{Nij3} that $N$ should satisfy the following condition
\begin{equation}
	\label{eq:Nijenhuis condition1} N\big(N(a)\circ_\lambda b+a\circ_\lambda N(b)-N(a\circ_\lambda b)\big)=N (a)\circ_\lambda N(b).
\end{equation}

It follows from \eqref{Nij4} and \eqref{Nij5} that $N$ should satisfy the following condition
\begin{equation}
	\label{eq:Nijenhuis condition2} N\big([N(a)_{\lambda}b]+[a_{\lambda}N(b)]-N[a_{\lambda}b]\big)=[N(a)_{\lambda}N(b)].
\end{equation}

Recall from \cite{Yuan} that a {\bf Nijenhuis operator} on an associative conformal algebra $(P,\circ_\lambda)$ is a $\Comp[\partial]$-module homomorphism $N:P\rightarrow P$ satisfying \eqref{eq:Nijenhuis condition1}.

Recall from  \cite{LZY} that a {\bf Nijenhuis operator} on a Lie conformal algebra $(P,[\cdot_\lambda\cdot])$ is a $\Comp[\partial]$-module homomorphism $N:P\rightarrow P$ satisfying \eqref{eq:Nijenhuis condition2}.

\begin{defi}
Let $(P,\circ_{\lambda},[\cdot_{\lambda}\cdot])$ be a noncommutative Poisson conformal algebra. A $\Comp[\partial]$-module homomorphism $N:P\rightarrow P$ is called a {\bf Nijenhuis operator} on $P$ if $N$ is both a Nijenhuis operator  on the associative conformal algebra $(P,\circ_\lambda)$ and a Nijenhuis operator  on the  Lie conformal algebra $(P,[\cdot_\lambda\cdot])$.
\end{defi}	

\begin{pro}\label{eq:Nijenhuis properties}
Let $(P,\circ_{\lambda},[\cdot_{\lambda}\cdot])$ be a noncommutative Poisson conformal algebra. If $N:P\rightarrow P$ is a Nijenhuis operator on the noncommutative Poisson conformal algebra $P$, then $(P,\circ_{N\lambda},[\cdot_{\lambda}\cdot]_N)$ is a noncommutative Poisson conformal algebra, where $\circ_{N\lambda}$ and $[\cdot_\lambda\cdot]_N$ are  defined by
\begin{eqnarray}
\label{eq:deformed operation}	a\circ_{N\lambda}b&=&N(a)\circ_\lambda b+a\circ_\lambda N(b)-N(a\circ_\lambda b),\\
\label{eq:deformed bracket}	[a_{\lambda}b]_{N}&=&[N(a)_{\lambda}b]+[a_{\lambda}N(b)]-N[a_{\lambda}b],\quad\forall~a,b\in P.
\end{eqnarray}
Furthermore, $N$ is a noncommutative Poisson conformal algebra homomorphism from $(P,\circ_{N\lambda},[\cdot_{\lambda}\cdot]_N)$ to $(P,\circ_{\lambda},[\cdot_{\lambda}\cdot])$.
\end{pro}
\begin{proof}
Since $N$ is a Nijenhuis operator on the associative conformal algebra $(P,\circ_\lambda)$, $(P,\circ_{N\lambda})$ is an associative conformal algebra (\cite{Yuan}).  Since $N$ is a Nijenhuis operator on the Lie conformal algebra $(P,[\cdot_\lambda\cdot])$, $(P,[\cdot_\lambda\cdot]_N)$ is a Lie conformal algebra  (\cite{LZY}).  It is straightforward to check that
$$	[a_{\lambda}(b\circ_{N\mu}c)]_N=[a_{\lambda}b]_N\circ_{N\lambda+\mu}c+b\circ _{N\mu}[a_{\lambda}c]_N.$$
Thus $(P,\circ_{N\lambda},[\cdot_{\lambda}\cdot]_N)$ is a noncommutative Poisson conformal algebra. The rest is direct.
\end{proof}

A trivial linear deformation of a noncommutative Poisson conformal algebra $(P,\circ_{\lambda},[\cdot_{\lambda}\cdot])$  gives rise to a Nijenhuis operator on the noncommutative Poisson conformal algebra $P$. In fact, the converse is also true.

\begin{thm}
	Let $(P,\circ_{\lambda},[\cdot_{\lambda}\cdot])$ be a noncommutative Poisson conformal algebra. If $N:P\rightarrow P$ is a Nijenhuis operator on the noncommutative Poisson conformal algebra $P$, then a linear deformation of the noncommutative Poisson conformal algebra $P$ can be obtained by putting
	\begin{eqnarray}
	\{a,b\}_\varpi&=&N(a)\circ_\lambda b+a\circ_\lambda N(b)-N(a\circ_\lambda b),\\
	\{a,b\}_\omega&=&[N(a)_{\lambda}b]+[a_{\lambda}N(b)]-N[a_{\lambda}b],\quad\forall~a,b\in P.
	\end{eqnarray}	
Furthermore, this linear deformation is trivial.
	\end{thm}
\begin{proof}
	Since $N$ is a Nijenhuis operator on the noncommutative Poisson conformal algebra $P$, by Proposition \ref{eq:Nijenhuis properties}, $(P,\{\cdot_\lambda\cdot\}_\varpi,\{\cdot_\lambda\cdot\}_\omega)$ is a noncommutative Poisson conformal algebra, which implies that  \eqref{eq:ass conformal1}-\eqref{eq:Poisson Leibniz}  hold.
	
	Since $\omega+\varpi=\dM_{\FGV} N$, $\omega+\varpi$ is a $2$-cocycle and thus \eqref{eq:deformation-FGV1}-\eqref{eq:deformation-FGV3} hold.
	
	It is not hard to see that \eqref{Nij1}-\eqref{Nij5} hold. Thus this linear deformation is trivial.
\end{proof}

\section*{Appendix}

{\bf The proof of Lemma \ref{lem:FGV2}}:

First, we show that $\dM^{m,1}_{\rm H}\circ \dM^{m,0}_{\rm CE}=\dM^{m-1,2}_{\rm CE}\circ \dM^{m-1,1}_{\rm H}$ for $m\geq 1$. By a direct calculation, we have
\begin{eqnarray*}
	&&\{x_{1\lambda_{1}}\cdots x_{m \lambda_{m}}a_{1\mu_1}a_{2}\}_{\dM^{m,1}_{\rm H}\circ \dM^{m,0}_{\rm CE} \gamma}\\
	%&=&{a_{1}}\circ_{\mu_{1}}\{x_{1\lambda_{1}}\cdots x_{m \lambda_{m}}a_{2}\}_{\dM^{m,0}_{\rm CE} \gamma}
	%-\{{x_{1\lambda_{1}}\cdots x_{m \lambda_{m}} ({a_{1}}\circ_{\mu_{1}} a_{2})}\}_{\dM^{m,0}_{\rm CE} \gamma}
	%+{\{x_{1\lambda_{1}}\cdots x_{m \lambda_{m}}{a_{1}}\}_{\dM^{m,0}_{\rm CE} %\gamma}}\circ_{\sum_{k=1}^{m}\lambda_{k}+\mu_{1}} a_{2}\\
	&=&\sum_{i=1}^{m}(-1)^{i+1}a_{1}\circ_{\mu_{1}}(x_{i\lambda_{i}}\{x_{1\lambda_{1}}\cdots\hat{x_{i\lambda_{i}}}\cdots x_{m\lambda_{m}}a_{2}\}_{\gamma})\\ &&+\sum_{i,j=1,i<j}^{m}(-1)^{m+i+j+1}a_{1}\circ_{\mu_{1}}\{x_{1\lambda_{1}}\cdots\hat{x_{i\lambda_{i}}}\cdots\hat{x_{j\lambda_{j}}}\cdots x_{m\lambda_{m}}a_{2{-\sum_{k=1}^{m}\lambda_{k}-\partial}}[x_{i\lambda_{i}}x_{j}]\}_{\gamma}\\      &&+(-1)^{m}a_{1}\circ_{\mu_{1}}(a_{2{-\sum_{k=1}^{m}\lambda_{k}-\partial}}\{x_{1\lambda_{1}}\cdots x_{m}\}_{\gamma})\\
	&&+\sum_{i=1}^{m}(-1)^{i}a_{1}\circ_{\mu_{1}}\{x_{1\lambda_{1}}\cdots\hat{x_{i}}\cdots x_{m\lambda_{m}}[x_{i\lambda_{i}}a_{2}]\}_{\gamma}\\
	&&-\sum_{i=1}^{m}(-1)^{i+1}x_{i\lambda_{i}}\{x_{1\lambda_{1}}\cdots\hat{x_{i\lambda_{i}}}\cdots x_{m\lambda_{m}}({a_{1}}\circ_{\mu_{1}} a_{2})\}_{\gamma}\\ &&-\sum_{i,j=1,i<j}^{m}(-1)^{m+i+j+1}\{x_{1\lambda_{1}}\cdots\hat{x_{i\lambda_{i}}}\cdots\hat{x_{j\lambda_{j}}}\cdots x_{m\lambda_{m}}({a_{1}}\circ_{\mu_{1}} a_{2})_{-\sum_{k=1}^{m}\lambda_{k}-\partial}[x_{i\lambda_{i}}x_{j}]\}_{\gamma}\\&&-(-1)^{m}({a_{1}}\circ_{\mu_{1}} a_{2})_{-\sum_{k=1}^{m}\lambda_{k}-\partial}\{x_{1\lambda_{1}}\cdots x_{m}\}_{\gamma}\\
	&&-\sum_{i=1}^{m}(-1)^{i}\{x_{1\lambda_{1}}\cdots\hat{x_{i}}\cdots x_{m\lambda_{m}}[x_{i\lambda_{i}}({a_{1}}\circ_{\mu_{1}} a_{2})]\}_{\gamma}\\
	&&+\sum_{i=1}^{m}(-1)^{i+1}x_{i\lambda_{i}}\{x_{1\lambda_{1}}\cdots\hat{x_{i\lambda_{i}}}\cdots x_{m\lambda_{m}}a_{1}\}_{\gamma}\circ_{\sum_{k=1}^{m}\lambda_{k}+\mu_{1}}a_{2}\\ &&+\sum_{i,j=1,i<j}^{m}(-1)^{m+i+j+1}\{x_{1\lambda_{1}}\cdots\hat{x_{i\lambda_{i}}}\cdots\hat{x_{j\lambda_{j}}}\cdots x_{m\lambda_{m}}a_{1{-\sum_{l=1}^{m}\lambda_{l}-\partial}}[x_{i\lambda_{i}}x_{j}]\}_{\gamma}\circ_{\sum_{k=1}^{m}\lambda_{k}+\mu_{1}}a_{2}\\      &&+(-1)^{m}(a_{1{-\sum_{l=1}^{m}\lambda_{l}-\partial}}\{x_{1\lambda_{1}}\cdots x_{m}\}_{\gamma})\circ_{\sum_{k=1}^{m}\lambda_{k}+\mu_{1}}a_{2}\\
	&&+\sum_{i=1}^{m}(-1)^{i}\{x_{1\lambda_{1}}\cdots\hat{x_{i}}\cdots x_{m\lambda_{m}}[x_{i\lambda_{i}}a_{1}]\}_{\gamma}\circ_{\sum_{k=1}^{m}\lambda_{k}+\mu_{1}}a_{2},
\end{eqnarray*}
and
\begin{eqnarray*}
	&&\{x_{1\lambda_{1}}\cdots x_{m \lambda_{m}}a_{1\mu_1}a_{2}\}_{\dM^{m-1,2}_{\rm CE}\circ \dM^{m-1,1}_{\rm H} \gamma}\\
	%&=&\sum_{i=1}^{m}(-1)^{i+1}\Big(x_{i\lambda_{i}}\{x_{1\lambda_{1}}\cdots\hat{x_{i\lambda_{i}}}\cdots x_{m \lambda_{m}}a_{1\mu_1}a_{2}\}_{\dM^{m-1,1}_{\rm H} \gamma}
	%-\{x_{1\lambda_{1}}\cdots\hat{x_{i\lambda_{i}}}\cdots x_{m \lambda_{m}}([x_{i\lambda_i} a_1])_{\lambda_i+\mu_1}a_{2}\}_{\dM^{m-1,1}_{\rm H} \gamma}\\
	%&&-\{x_{1\lambda_{1}}\cdots\hat{x_{i\lambda_{i}}}\cdots x_{m \lambda_{m}}a_{1\mu_{1}}([x_{i\lambda_i} a_2])\}_{\dM^{m-1,1}_{\rm H} \gamma}\Big)
	%	+\sum_{i,j=1,i<j}^{m}(-1)^{i+j}\{([x_{i\lambda_i}x_j])_{\lambda_i+\lambda_j}x_{1\lambda_{1}}\cdots\hat{x_{i\lambda_i}}\cdots\hat{x_{j\lambda_j}}\cdots x_{m \lambda_{m}}a_{1\mu_1}a_{2}\}_{\dM^{m-1,1}_{\rm H} \gamma}\\
	&=&\sum_{i=1}^{m}(-1)^{i+1}\Big(x_{i\lambda_{i}}({a_{1}}\circ_{\mu_{1}}\{x_{1\lambda_{1}}\cdots \hat{x_{i\lambda_{i}}}\cdots x_{m \lambda_{m}}a_{2}\}_{\gamma})
	-x_{i\lambda_{i}}\{x_{1\lambda_{1}}\cdots\hat{x_{i\lambda_{i}}}\cdots x_{m \lambda_{m}}({a_{1}}\circ_{\mu_{1}} a_{2})\}_{\gamma}\\
	&&+x_{i\lambda_{i}}(\{x_{1\lambda_{1}}\cdots \hat{x_{i\lambda_{i}}}\cdots x_{m \lambda_{m}}a_{1}\}_{\gamma}\circ_{\sum_{k=1}^{m}\lambda_{k}-\lambda_{i}+\mu_{1}}a_{2})\Big)\\
	&&-\sum_{i=1}^{m}(-1)^{i+1}\Big(([x_{i\lambda_i} a_1])\circ_{\lambda_i+\mu_1}\{x_{1\lambda_{1}}\cdots \hat{x_{i\lambda_{i}}}\cdots x_{m \lambda_{m}}a_{2}\}_{\gamma}
	-\{x_{1\lambda_{1}}\cdots\hat{x_{i\lambda_{i}}}\cdots x_{m \lambda_{m}}(([x_{i\lambda_i} a_1])\circ_{\lambda_i+\mu_1} a_{2})\}_{\gamma}\\
	&&+\{x_{1\lambda_{1}}\cdots \hat{x_{i\lambda_{i}}}\cdots x_{m \lambda_{m}}([x_{i\lambda_i} a_1])\}_{\gamma}\circ_{\sum_{k=1}^{m}\lambda_{k}+\mu_{1}}a_{2}\Big)\\
	&&-\sum_{i=1}^{m}(-1)^{i+1}\Big({a_{1}}\circ_{\mu_{1}}\{x_{1\lambda_{1}}\cdots \hat{x_{i\lambda_{i}}}\cdots x_{m \lambda_{m}}([x_{i\lambda_i} a_2])\}_{\gamma}
	-\{x_{1\lambda_{1}}\cdots\hat{x_{i\lambda_{i}}}\cdots x_{m \lambda_{m}}({a_{1}}\circ_{\mu_{1}} [x_{i\lambda_i} a_2])\}_{\gamma}\\
	&&+\{x_{1\lambda_{1}}\cdots \hat{x_{i\lambda_{i}}}\cdots x_{m \lambda_{m}}a_{1}\}_{\gamma}\circ_{\sum_{k=1}^{m}\lambda_{k}-\lambda_{i}+\mu_{1}}([x_{i\lambda_i} a_2])\Big)\\
	&&+\sum_{i,j=1,i<j}^{m}(-1)^{i+j}\Big(a_{1}\circ_{\mu_{1}}\{ ([x_{i\lambda_i}x_j])_{\lambda_i+\lambda_j}x_{1\lambda_{1}}\cdots\hat{x_{i\lambda_i}}\cdots\hat{x_{j\lambda_j}}\cdots x_{m \lambda_{m}}a_{2}\}_{\gamma}\\
	&&-\{([x_{i\lambda_i}x_j])_{\lambda_i+\lambda_j}x_{1\lambda_{1}}\cdots\hat{x_{i\lambda_i}}\cdots\hat{x_{j\lambda_j}}\cdots x_{m \lambda_{m}}(a_{1}\circ_{\mu_{1}}a_{2})\}_{\gamma}\\
	&&+\{([x_{i\lambda_i}x_j])_{\lambda_i+\lambda_j}x_{1\lambda_{1}}\cdots\hat{x_{i\lambda_i}}\cdots\hat{x_{j\lambda_j}}\cdots x_{m \lambda_{m}}a_{1}\}_{\gamma}\circ_{\sum_{k=1}^{m}\lambda_{k}+\mu_{1}}a_{2}\Big).
\end{eqnarray*}
Then, by \eqref{eq:Poisson conformal Leibniz}, \eqref{eq:Poisson module 1}, \eqref{eq:Poisson module 2} and  \eqref{eq:Poisson module 3}, we have
\begin{eqnarray*}
	&&\{x_{1\lambda_{1}}\cdots x_{m \lambda_{m}}a_{1\mu_1}a_{2}\}_{\dM^{m,1}_{\rm H}\circ \dM^{m,0}_{\rm CE} \gamma}
	-\{x_{1\lambda_{1}}\cdots x_{m \lambda_{m}}a_{1\mu_1}a_{2}\}_{\dM^{m-1,2}_{\rm CE}\circ \dM^{m-1,1}_{\rm H} \gamma}\\
	&=&\sum_{i=1}^{m}(-1)^{i+1}\Big(a_{1}\circ_{\mu_{1}}(x_{i\lambda_{i}}\{x_{1\lambda_{1}}\cdots\hat{x_{i\lambda_{i}}}\cdots x_{m\lambda_{m}}a_{2}\}_{\gamma})
	-x_{i\lambda_{i}}({a_{1}}\circ_{\mu_{1}}\{x_{1\lambda_{1}}\cdots \hat{x_{i\lambda_{i}}}\cdots x_{m \lambda_{m}}a_{2}\}_{\gamma})\\
	&&+([x_{i\lambda_i} a_1])\circ_{\lambda_i+\mu_1}\{x_{1\lambda_{1}}\cdots \hat{x_{i\lambda_{i}}}\cdots x_{m \lambda_{m}}a_{2}\}_{\gamma}\Big)\\
	&&+\Big(\sum_{i,j=1,i<j}^{m}(-1)^{m+i+j+1}a_{1}\circ_{\mu_{1}}\{x_{1\lambda_{1}}\cdots\hat{x_{i\lambda_{i}}}\cdots\hat{x_{j\lambda_{j}}}\cdots x_{m\lambda_{m}}a_{2{-\sum_{k=1}^{m}\lambda_{k}-\partial}}[x_{i\lambda_{i}}x_{j}]\}_{\gamma}\\
	&&-\sum_{i,j=1,i<j}^{m}(-1)^{i+j}a_{1}\circ_{\mu_{1}}\{([x_{i\lambda_i}x_j])_{\lambda_i+\lambda_j}x_{1\lambda_{1}}\cdots\hat{x_{i\lambda_i}}\cdots\hat{x_{j\lambda_j}}\cdots x_{m \lambda_{m}}a_{2}\}_{\gamma}\Big)\\
	&&+(-1)^{m}\Big(a_{1}\circ_{\mu_{1}}(a_{2{-\sum_{k=1}^{m}\lambda_{k}-\partial}}\{x_{1\lambda_{1}}\cdots x_{m}\}_{\gamma})
	-({a_{1}}\circ_{\mu_{1}} a_{2})_{-\sum_{k=1}^{m}\lambda_{k}-\partial}\{x_{1\lambda_{1}}\cdots x_{m}\}_{\gamma}\\
	&&+(a_{1{-\sum_{l=1}^{m}\lambda_{l}-\partial}}\{x_{1\lambda_{1}}\cdots x_{m}\}_{\gamma})\circ_{\sum_{k=1}^{m}\lambda_{k}+\mu_{1}}a_{2}\Big)\\
	&&-\Big(\sum_{i,j=1,i<j}^{m}(-1)^{m+i+j+1}\{x_{1\lambda_{1}}\cdots\hat{x_{i\lambda_{i}}}\cdots\hat{x_{j\lambda_{j}}}\cdots x_{m\lambda_{m}}({a_{1}}\circ_{\mu_{1}} a_{2})_{-\sum_{k=1}^{m}\lambda_{k}-\partial}[x_{i\lambda_{i}}x_{j}]\}_{\gamma}\\
	&&-\sum_{i,j=1,i<j}^{m}(-1)^{i+j}\{([x_{i\lambda_i}x_j])_{\lambda_i+\lambda_j}x_{1\lambda_{1}}\cdots\hat{x_{i\lambda_i}}\cdots\hat{x_{j\lambda_j}}\cdots x_{m \lambda_{m}}(a_{1}\circ_{\mu_{1}}a_{2})\}_{\gamma}\Big)\\
	&&+\sum_{i=1}^{m}(-1)^{i+1}\Big(\{x_{1\lambda_{1}}\cdots\hat{x_{i}}\cdots x_{m\lambda_{m}}[x_{i\lambda_{i}}({a_{1}}\circ_{\mu_{1}} a_{2})]\}_{\gamma}
	-\{x_{1\lambda_{1}}\cdots\hat{x_{i\lambda_{i}}}\cdots x_{m \lambda_{m}}(([x_{i\lambda_i} a_1])\circ_{\lambda_i+\mu_1} a_{2})\}_{\gamma}\\
	&&-\{x_{1\lambda_{1}}\cdots\hat{x_{i\lambda_{i}}}\cdots x_{m \lambda_{m}}({a_{1}}\circ_{\mu_{1}} [x_{i\lambda_i} a_2])\}_{\gamma}\Big)\\
	&&+\sum_{i=1}^{m}(-1)^{i+1}\Big(x_{i\lambda_{i}}\{x_{1\lambda_{1}}\cdots\hat{x_{i\lambda_{i}}}\cdots x_{m\lambda_{m}}a_{1}\}_{\gamma}\circ_{\sum_{k=1}^{m}\lambda_{k}+\mu_{1}}a_{2}
	-x_{i\lambda_{i}}(\{x_{1\lambda_{1}}\cdots \hat{x_{i\lambda_{i}}}\cdots x_{m \lambda_{m}}a_{1}\}_{\gamma}\circ_{\sum_{k=1}^{m}\lambda_{k}-\lambda_{i}+\mu_{1}}a_{2})\\
	&&+\{x_{1\lambda_{1}}\cdots \hat{x_{i\lambda_{i}}}\cdots x_{m \lambda_{m}}a_{1}\}_{\gamma}\circ_{\sum_{k=1}^{m}\lambda_{k}-\lambda_{i}+\mu_{1}}([x_{i\lambda_i} a_2])\Big)\\
	&&+\Big(\sum_{i,j=1,i<j}^{m}(-1)^{m+i+j+1}\{x_{1\lambda_{1}}\cdots\hat{x_{i\lambda_{i}}}\cdots\hat{x_{j\lambda_{j}}}\cdots x_{m\lambda_{m}}a_{1{-\sum_{l=1}^{m}\lambda_{l}-\partial}}[x_{i\lambda_{i}}x_{j}]\}_{\gamma}\circ_{\sum_{k=1}^{m}\lambda_{k}+\mu_{1}}a_{2}\\
	&&-\sum_{i,j=1,i<j}^{m}\{([x_{i\lambda_i}x_j])_{\lambda_i+\lambda_j}x_{1\lambda_{1}}\cdots\hat{x_{i\lambda_i}}\cdots\hat{x_{j\lambda_j}}\cdots x_{m \lambda_{m}}a_{1}\}_{\gamma}\circ_{\sum_{k=1}^{m}\lambda_{k}+\mu_{1}}a_{2}\Big)\\
	&=&0,
\end{eqnarray*}
which implies that $ \dM^{m,1}_{\rm H}\circ \dM^{m,0}_{\rm CE}=\dM^{m-1,2}_{\rm CE}\circ \dM^{m-1,1}_{\rm H} $ for $m\geq 1$.

Next, we show that $\dM^{m+1,n}_{\rm H}\circ \dM^{m,n}_{\rm CE}=\dM^{m,n+1}_{\rm CE}\circ \dM^{m,n}_{\rm H}$ for $m\geq 1, n\geq 2$. By a direct calculation, we have
\begin{eqnarray*}
	&&\{x_{1\lambda_{1}}\cdots x_{m+1 \lambda_{m+1}}a_{1\mu_1}\cdots a_{n\mu_{n}}a_{n+1}\}_{\dM^{m+1,n}_{\rm H}\circ \dM^{m,n}_{\rm CE} \gamma}\\
	%&=&{a_{1}}\circ_{\mu_{1}}\{x_{1\lambda_{1}}\cdots x_{m+1 \lambda_{m+1}}a_{2\mu_2}\cdots a_{n\mu_{n}}a_{n+1}\}_{\dM^{m,n}_{\rm CE} \gamma}\\
	%&&+\sum_{i=1}^{n}(-1)^{i} \{{x_{1\lambda_{1}}\cdots x_{m+1 \lambda_{m+1}}a_{1}}_{\mu_{1}}\dots {a_{i-1}}_{\mu_{i-1}}({a_{i}}\circ_{\mu_{i}} a_{i+1})_{\mu_{i}+\mu_{i+1}}{a_{i+2}}_{\mu_{i+2}}\dots {a_{n}}_{\mu_{n}}{a_{n+1}}\}_{\dM^{m,n}_{\rm CE} \gamma}\\
	%&&+(-1)^{n+1}{\{x_{1\lambda_{1}}\cdots x_{m+1 \lambda_{m+1}}{a_{1}}_{\mu_{1}}\dots {a_{n-1}}_{\mu_{n-1}}{ a_{n}}\}_{\dM^{m,n}_{\rm CE} \gamma}}\circ_{\mu^\flat_{n+1}} a_{n+1}\\
	&=&\sum_{i=1}^{m+1}(-1)^{i+1}{a_{1}}\circ_{\mu_{1}}(x_{i \lambda_{i}}\{x_{1\lambda_{1}}\cdots \hat{x_{i\lambda_{i}}}\cdots x_{m+1 \lambda_{m+1}}a_{2\mu_2}\cdots a_{n\mu_{n}}a_{n+1}\}_{\gamma})\\
	&&-\sum_{i=1}^{m+1}(-1)^{i+1}\sum_{j=2}^{n+1}{a_{1}}\circ_{\mu_{1}}\{x_{1\lambda_{1}}\cdots\hat{x_{i\lambda_{i}}}\cdots x_{m+1 \lambda_{m+1}}a_{2\mu_{2}}\cdots([x_{i\lambda_i} a_j])_{\lambda_i+\mu_j}\cdots a_{n \mu_{n}}a_{n+1}\}_{\gamma}\\
	&&+\sum_{i,j=1,i<j}^{m+1}(-1)^{i+j}{a_{1}}\circ_{\mu_{1}}\{([x_{i\lambda_i}x_j])_{\lambda_i+\lambda_j}x_{1\lambda_{1}}\cdots\hat{x_{i\lambda_i}}\cdots\hat{x_{j\lambda_j}}\cdots x_{m+1 \lambda_{m+1}}a_{2\mu_2}\cdots a_{n\mu_{n}}a_{n+1}\}_{\gamma}\\
	&&+\sum_{i=1}^{n}(-1)^{i}\sum_{j=1}^{m+1}(-1)^{j+1}x_{j \lambda_{j}}\{x_{1\lambda_{1}}\cdots \hat{x_{j\lambda_{j}}}\cdots x_{m+1 \lambda_{m+1}}a_{1\mu_1}\cdots ({a_{i}}\circ_{\mu_{i}} a_{i+1})_{\mu_{i}+\mu_{i+1}}\cdots a_{n\mu_{n}}a_{n+1}\}_{\gamma}\\
	&&-\sum_{i=1}^{n}(-1)^{i}\sum_{j=1}^{m+1}(-1)^{j+1}\sum_{k=1}^{i-1}\{x_{1\lambda_{1}}\cdots\hat{x_{j\lambda_{j}}}\cdots x_{m+1 \lambda_{m+1}}a_{1\mu_{1}}\cdots([x_{j\lambda_j} a_k])_{\lambda_j+\mu_k}\cdots ({a_{i}}\circ_{\mu_{i}} a_{i+1})_{\mu_{i}+\mu_{i+1}}\cdots a_{n+1}\}_{\gamma}\\
	&&-\sum_{i=1}^{n}(-1)^{i}\sum_{j=1}^{m+1}(-1)^{j+1}\{x_{1\lambda_{1}}\cdots\hat{x_{j\lambda_{j}}}\cdots x_{m+1 \lambda_{m+1}}a_{1\mu_{1}}\cdots([x_{j\lambda_j} ({a_{i}}\circ_{\mu_{i}} a_{i+1})])_{\lambda_j+\mu_{i}+\mu_{i+1}}\cdots a_{n+1}\}_{\gamma}\\
	&&-\sum_{i=1}^{n}(-1)^{i}\sum_{j=1}^{m+1}(-1)^{j+1}\sum_{k=i+2}^{n+1}\{x_{1\lambda_{1}}\cdots\hat{x_{j\lambda_{j}}}\cdots x_{m+1 \lambda_{m+1}}a_{1\mu_{1}}\cdots({a_{i}}\circ_{\mu_{i}} a_{i+1})_{\mu_{i}+\mu_{i+1}}\cdots ([x_{j\lambda_j} a_k])_{\lambda_j+\mu_k}\cdots a_{n+1}\}_{\gamma}\\
	&&+\sum_{i=1}^{n}(-1)^{i}\sum_{j,k=1,j<k}^{m+1}(-1)^{j+k}\{([x_{j\lambda_j}x_k])_{\lambda_j+\lambda_k}x_{1\lambda_{1}}\cdots\hat{x_{j\lambda_j}}\cdots\hat{x_{k\lambda_k}}\cdots x_{m+1 \lambda_{m+1}}a_{1\mu_1}\cdots({a_{i}}\circ_{\mu_{i}} a_{i+1})_{\mu_{i}+\mu_{i+1}}\cdots a_{n+1}\}_{\gamma}\\
	&&+(-1)^{n+1}\sum_{i=1}^{m+1}(-1)^{i+1}x_{i \lambda_{i}}{\{x_{1\lambda_{1}}\cdots \hat{x_{i\lambda_{i}}}\cdots x_{m+1 \lambda_{m+1}}{a_{1}}_{\mu_{1}}\dots {a_{n-1}}_{\mu_{n-1}}{ a_{n}}\}_{ \gamma}}\circ_{\mu^\flat_{n+1}} a_{n+1}\\
	&&-(-1)^{n+1}\sum_{i=1}^{m+1}(-1)^{i+1}\sum_{j=1}^{n}{\{x_{1\lambda_{1}}\cdots \hat{x_{i\lambda_{i}}}\cdots x_{m+1 \lambda_{m+1}}{a_{1}}_{\mu_{1}}\dots ([x_{i \lambda_{i}}a_{j}])_{\lambda_{i}+\mu_{j}}\cdots {a_{n-1}}_{\mu_{n-1}}{ a_{n}}\}_{ \gamma}}\circ_{\mu^\flat_{n+1}} a_{n+1}\\
	&&+(-1)^{n+1}\sum_{i,j=1,i<j}^{m+1}(-1)^{i+j}{\{([x_{i \lambda_{i}}x_{j}])_{\lambda_{i}+\lambda_{j}}x_{1\lambda_{1}}\cdots \hat{x_{i\lambda_{i}}}\cdots \hat{x_{j\lambda_{j}}} \cdots x_{m+1 \lambda_{m+1}}{a_{1}}_{\mu_{1}}\dots {a_{n-1}}_{\mu_{n-1}}{ a_{n}}\}_{ \gamma}}\circ_{\mu^\flat_{n+1}} a_{n+1},
\end{eqnarray*}
and
\begin{eqnarray*}
	&&\{x_{1\lambda_{1}}\cdots x_{m+1 \lambda_{m+1}}a_{1\mu_1}\cdots a_{n\mu_{n}}a_{n+1}\}_{\dM^{m,n+1}_{\rm CE}\circ \dM^{m,n}_{\rm H} \gamma}\\
	%&=&\sum_{i=1}^{m+1}(-1)^{i+1}x_{i\lambda_{i}}\{x_{1\lambda_{1}}\cdots\hat{x_{i\lambda_{i}}}\cdots x_{m+1 \lambda_{m+1}}a_{1\mu_1}\cdots a_{n\mu_{n}}a_{n+1}\}_{\dM^{m,n}_{\rm H} \gamma}\\
	%&&-\sum_{i=1}^{m+1}(-1)^{i+1}\sum_{j=1}^{n+1}\{x_{1\lambda_{1}}\cdots\hat{x_{i\lambda_{i}}}\cdots x_{m+1 \lambda_{m+1}}a_{1\mu_{1}}\cdots([x_{i\lambda_i} a_j])_{\lambda_i+\mu_j}\cdots a_{n \lambda_{n}}a_{n+1}\}_{\dM^{m,n}_{\rm H} \gamma}\\
	%&&+\sum_{i,j=1,i<j}^{m+1}(-1)^{i+j}\{([x_{i\lambda_i}x_j])_{\lambda_i+\lambda_j}x_{1\lambda_{1}}\cdots\hat{x_{i\lambda_i}}\cdots\hat{x_{j\lambda_j}}\cdots x_{m+1 \lambda_{m+1}}a_{1\mu_1}\cdots a_{n\mu_{n}}a_{n+1}\}_{\dM^{m,n}_{\rm H} \gamma}\\
	&=&\sum_{i=1}^{m+1}(-1)^{i+1}x_{i \lambda_{i}}({a_{1}}\circ_{\mu_{1}}\{x_{1\lambda_{1}}\cdots \hat{x_{i\lambda_{i}}}\cdots x_{m+1 \lambda_{m+1}}a_{2\mu_2}\cdots a_{n\mu_{n}}a_{n+1}\}_{\gamma})\\
	&&+\sum_{i=1}^{m+1}(-1)^{i+1}\sum_{j=1}^{n}(-1)^{j}x_{i \lambda_{i}}\{x_{1\lambda_{1}}\cdots \hat{x_{i\lambda_{i}}}\cdots x_{m+1 \lambda_{m+1}}a_{1\mu_1}\cdots ({a_{j}}\circ_{\mu_{j}} a_{j+1})_{\mu_{j}+\mu_{j+1}}\cdots a_{n\mu_{n}}a_{n+1}\}_{\gamma}\\
	&&+(-1)^{n+1}\sum_{i=1}^{m+1}(-1)^{i+1}x_{i \lambda_{i}}{\{x_{1\lambda_{1}}\cdots \hat{x_{i\lambda_{i}}}\cdots x_{m+1 \lambda_{m+1}}{a_{1}}_{\mu_{1}}\dots {a_{n-1}}_{\mu_{n-1}}{ a_{n}}\}_{ \gamma}}\circ_{\mu^\flat_{n+1}-\lambda_{i}} a_{n+1}\\
	&&-\sum_{i=1}^{m+1}(-1)^{i+1}[x_{i \lambda_{i}}{a_{1}}]\circ_{\lambda_{i}+\mu_{1}}\{x_{1\lambda_{1}}\cdots \hat{x_{i\lambda_{i}}}\cdots x_{m+1 \lambda_{m+1}}a_{2\mu_2}\cdots a_{n\mu_{n}}a_{n+1}\}_{\gamma}\\
	&&-\sum_{i=1}^{m+1}(-1)^{i+1}\sum_{j=1}^{n}(-1)^{j}\{x_{1\lambda_{1}}\cdots\hat{x_{i\lambda_{i}}}\cdots x_{m+1 \lambda_{m+1}}a_{1\mu_{1}}\cdots([x_{i\lambda_i} a_{j}]\circ_{\lambda_{i}+\mu_{j}} a_{j+1})_{\lambda_{i}+\mu_{j}+\mu_{j+1}}\cdots a_{n+1}\}_{\gamma}\\
	&&-\sum_{i=1}^{m+1}(-1)^{i+1}\sum_{j=1}^{n}(-1)^{j}\{x_{1\lambda_{1}}\cdots\hat{x_{i\lambda_{i}}}\cdots x_{m+1 \lambda_{m+1}}a_{1\mu_{1}}\cdots(a_{j}\circ_{\mu_{j}}[x_{i\lambda_i} {a_{j+1}}])_{\lambda_{i}+\mu_{j}+\mu_{j+1}}\cdots a_{n+1}\}_{\gamma}\\
	&&-\sum_{i=1}^{m+1}(-1)^{i+1}\sum_{j=1}^{n}(-1)^{j}\sum_{k=1}^{j-1}\{x_{1\lambda_{1}}\cdots\hat{x_{i\lambda_{i}}}\cdots x_{m+1 \lambda_{m+1}}a_{1\mu_{1}}\cdots([x_{i\lambda_i} a_k])_{\lambda_i+\mu_k}\cdots ({a_{j}}\circ_{\mu_{j}} a_{j+1})_{\mu_{j}+\mu_{j+1}}\cdots a_{n+1}\}_{\gamma}\\
	&&-\sum_{i=1}^{m+1}(-1)^{i+1}\sum_{j=1}^{n}(-1)^{j}\sum_{k=j+2}^{n+1}\{x_{1\lambda_{1}}\cdots\hat{x_{i\lambda_{i}}}\cdots x_{m+1 \lambda_{m+1}}a_{1\mu_{1}}\cdots({a_{j}}\circ_{\mu_{j}} a_{j+1})_{\mu_{j}+\mu_{j+1}}\cdots ([x_{i\lambda_i} a_k])_{\lambda_i+\mu_k}\cdots a_{n+1}\}_{\gamma}\\
	&&-(-1)^{n+1}\sum_{i=1}^{m+1}(-1)^{i+1}\sum_{j=1}^{n}{\{x_{1\lambda_{1}}\cdots \hat{x_{i\lambda_{i}}}\cdots x_{m+1 \lambda_{m+1}}{a_{1}}_{\mu_{1}}\dots ([x_{i \lambda_{i}}a_{j}])_{\lambda_{i}+\mu_{j}}\cdots {a_{n-1}}_{\mu_{n-1}}{ a_{n}}\}_{ \gamma}}\circ_{\mu^\flat_{n+1}} a_{n+1}\\
	&&-\sum_{i=1}^{m+1}(-1)^{i+1}\sum_{j=2}^{n+1}{a_{1}}\circ_{\mu_{1}}\{x_{1\lambda_{1}}\cdots\hat{x_{i\lambda_{i}}}\cdots x_{m+1 \lambda_{m+1}}a_{2\mu_{2}}\cdots([x_{i\lambda_i} a_j])_{\lambda_i+\mu_j}\cdots a_{n \mu_{n}}a_{n+1}\}_{\gamma}\\
	&&-(-1)^{n+1}\sum_{i=1}^{m+1}(-1)^{i+1}{\{x_{1\lambda_{1}}\cdots \hat{x_{i\lambda_{i}}}\cdots x_{m+1 \lambda_{m+1}}{a_{1}}_{\mu_{1}}\cdots {a_{n-1}}_{\mu_{n-1}}{ a_{n}}\}_{ \gamma}}\circ_{\mu^\flat_{n+1}-\lambda_{i}} [x_{i \lambda_{i}}a_{n+1}]\\
	&&+\sum_{i,j=1,i<j}^{m+1}(-1)^{i+j}{a_{1}}\circ_{\mu_{1}}\{([x_{i\lambda_i}x_j])_{\lambda_i+\lambda_j}x_{1\lambda_{1}}\cdots\hat{x_{i\lambda_i}}\cdots\hat{x_{j\lambda_j}}\cdots x_{m+1 \lambda_{m+1}}a_{2\mu_2}\cdots a_{n\mu_{n}}a_{n+1}\}_{\gamma}\\
	&&+\sum_{k=1}^{n}(-1)^{k}\sum_{i,j=1,i<j}^{m+1}(-1)^{j+k}\{([x_{i\lambda_i}x_j])_{\lambda_i+\lambda_j}x_{1\lambda_{1}}\cdots\hat{x_{i\lambda_i}}\cdots\hat{x_{j\lambda_j}}\cdots x_{m+1 \lambda_{m+1}}a_{1\mu_1}\cdots({a_{k}}\circ_{\mu_{k}} a_{k+1})_{\mu_{k}+\mu_{k+1}}\cdots a_{n+1}\}_{\gamma}\\
	&&+(-1)^{n+1}\sum_{i,j=1,i<j}^{m+1}(-1)^{i+j}{\{([x_{i \lambda_{i}}x_{j}])_{\lambda_{i}+\lambda_{j}}x_{1\lambda_{1}}\cdots \hat{x_{i\lambda_{i}}}\cdots \hat{x_{j\lambda_{j}}} \cdots x_{m+1 \lambda_{m+1}}{a_{1}}_{\mu_{1}}\dots {a_{n-1}}_{\mu_{n-1}}{ a_{n}}\}_{ \gamma}}\circ_{\mu^\flat_{n+1}} a_{n+1}
\end{eqnarray*}
where $\mu^\flat_{n+1}=\sum_{i=1}^{m+1}\lambda_i+\sum_{j=1}^n\mu_j$.

Then, by \eqref{eq:Poisson conformal Leibniz}, \eqref{eq:Poisson module 1} and \eqref{eq:Poisson module 2}, we have
\begin{eqnarray*}
	&&\{x_{1\lambda_{1}}\cdots x_{m+1 \lambda_{m+1}}a_{1\mu_1}\cdots a_{n\mu_{n}}a_{n+1}\}_{\dM^{m+1,n}_{\rm H}\circ \dM^{m,n}_{\rm CE} \gamma}
	-\{x_{1\lambda_{1}}\cdots x_{m+1 \lambda_{m+1}}a_{1\mu_1}\cdots a_{n\mu_{n}}a_{n+1}\}_{\dM^{m,n+1}_{\rm CE}\circ \dM^{m,n}_{\rm H} \gamma}\\
	&=&\sum_{i=1}^{m+1}(-1)^{i+1}{a_{1}}\circ_{\mu_{1}}(x_{i \lambda_{i}}\{x_{1\lambda_{1}}\cdots \hat{x_{i\lambda_{i}}}\cdots x_{m+1 \lambda_{m+1}}a_{2\mu_2}\cdots a_{n\mu_{n}}a_{n+1}\}_{\gamma})\\
	&&-\sum_{i=1}^{m+1}(-1)^{i+1}x_{i \lambda_{i}}({a_{1}}\circ_{\mu_{1}}\{x_{1\lambda_{1}}\cdots \hat{x_{i\lambda_{i}}}\cdots x_{m+1 \lambda_{m+1}}a_{2\mu_2}\cdots a_{n\mu_{n}}a_{n+1}\}_{\gamma})\\
	&&+\sum_{i=1}^{m+1}(-1)^{i+1}[x_{i \lambda_{i}}{a_{1}}]\circ_{\lambda_{i}+\mu_{1}}\{x_{1\lambda_{1}}\cdots \hat{x_{i\lambda_{i}}}\cdots x_{m+1 \lambda_{m+1}}a_{2\mu_2}\cdots a_{n\mu_{n}}a_{n+1}\}_{\gamma}\\
	&&-\sum_{i=1}^{n}(-1)^{i}\sum_{j=1}^{m+1}(-1)^{j+1}\{x_{1\lambda_{1}}\cdots\hat{x_{j\lambda_{j}}}\cdots x_{m+1 \lambda_{m+1}}a_{1\mu_{1}}\cdots([x_{j\lambda_j} ({a_{i}}\circ_{\mu_{i}} a_{i+1})])_{\lambda_j+\mu_{i}+\mu_{i+1}}\cdots a_{n+1}\}_{\gamma}\\
	&&+\sum_{i=1}^{m+1}(-1)^{i+1}\sum_{j=1}^{n}(-1)^{j}\{x_{1\lambda_{1}}\cdots\hat{x_{i\lambda_{i}}}\cdots x_{m+1 \lambda_{m+1}}a_{1\mu_{1}}\cdots([x_{i\lambda_i} a_{j}]\circ_{\lambda_{i}+\mu_{j}} a_{j+1})_{\lambda_{i}+\mu_{j}+\mu_{j+1}}\cdots a_{n+1}\}_{\gamma}\\
	&&+\sum_{i=1}^{m+1}(-1)^{i+1}\sum_{j=1}^{n}(-1)^{j}\{x_{1\lambda_{1}}\cdots\hat{x_{i\lambda_{i}}}\cdots x_{m+1 \lambda_{m+1}}a_{1\mu_{1}}\cdots(a_{j}\circ_{\mu_{j}}[x_{i\lambda_i} {a_{j+1}}])_{\lambda_{i}+\mu_{j}+\mu_{j+1}}\cdots a_{n+1}\}_{\gamma}\\
	&&+(-1)^{n+1}\sum_{i=1}^{m+1}(-1)^{i+1}x_{i \lambda_{i}}{\{x_{1\lambda_{1}}\cdots \hat{x_{i\lambda_{i}}}\cdots x_{m+1 \lambda_{m+1}}{a_{1}}_{\mu_{1}}\dots {a_{n-1}}_{\mu_{n-1}}{ a_{n}}\}_{ \gamma}}\circ_{\mu^\flat_{n+1}} a_{n+1}\\
	&&-(-1)^{n+1}\sum_{i=1}^{m+1}(-1)^{i+1}x_{i \lambda_{i}}{\{x_{1\lambda_{1}}\cdots \hat{x_{i\lambda_{i}}}\cdots x_{m+1 \lambda_{m+1}}{a_{1}}_{\mu_{1}}\dots {a_{n-1}}_{\mu_{n-1}}{ a_{n}}\}_{ \gamma}}\circ_{\mu^\flat_{n+1}-\lambda_{i}} a_{n+1}\\
	&&+(-1)^{n+1}\sum_{i=1}^{m+1}(-1)^{i+1}{\{x_{1\lambda_{1}}\cdots \hat{x_{i\lambda_{i}}}\cdots x_{m+1 \lambda_{m+1}}{a_{1}}_{\mu_{1}}\cdots {a_{n-1}}_{\mu_{n-1}}{ a_{n}}\}_{ \gamma}}\circ_{\mu^\flat_{n+1}-\lambda_{i}} [x_{i \lambda_{i}}a_{n+1}]\\
	&=&0,
\end{eqnarray*}
which implies that $ \dM^{m+1,n}_{\rm H}\circ \dM^{m,n}_{\rm CE}=\dM^{m,n+1}_{\rm CE}\circ \dM^{m,n}_{\rm H} $ for $m\geq 1, n\geq 2$.

\vspace{3mm}

{\bf Acknowledgements.} We give our warmest thanks to P. S. Kolesnikov,  Yunhe Sheng and Chongying Dong for very useful comments and discussions. This research was supported by the National Key Research and Development Program of China (2021YFA1002000), the National Natural Science Foundation of China (11901501), the China Postdoctoral Science Foundation (2021M700750, 2022T150109) and the Fundamental Research Funds for the Central Universities (2412022QD033).


\vspace{3mm}
\noindent



\begin{thebibliography}{abc}

\bibitem{Arn78}V.I. Arnol'd, Mathematical Methods of Classical Mechanics. \emph{Graduate Texts  Math.} 60, Springer-Verlag, New York-Heidelberg, 1978.


	
\bibitem{BAK}
	B. Bakalov, A. D'Andrea and V. G. Kac, Theory of finite pseudoalgebras. \emph{ Adv. Math.} 162 (2001),  1-140.
	
	
\bibitem{BDK}
B. Bakalov, A. De Sole and V. G. Kac, Computation of cohomology of Lie conformal and Poisson vertex algebras. \emph{Selecta Math. (N. S.)} 26 (2020), Paper No. 50.

\bibitem{BDK2}
B. Bakalov, A. De Sole and V. G. Kac, Computation of cohomology of vertex algebras. \emph{Jpn. J. Math.} 16 (2021), 81-154.

\bibitem{Bakalov}
B. Bakalov and V. G. Kac, Field algebras. \emph{Int. Math. Res. Not.} (2003), 123-159.

\bibitem{BKV}
B. Bakalov, V. G. Kac and A. A. Voronov, Cohomology of conformal algebras. \emph{Comm. Math. Phys.} 200 (1999), 561-598.

\bibitem{CohomologyPA1}
Y. Bao and Y. Ye, Cohomology structure for a Poisson algebra: I.  \emph{J. Algebra Appl.} 15 (2016), 1650034.

\bibitem{CohomologyPA2}
Y. Bao and Y. Ye,  Cohomology structure for a Poisson algebra: II. \emph{Sci. China Math.} 64 (2021), 903-920.

\bibitem{BSK}
A. Barakat, A. De Sole and V. G. Kac, Poisson vertex algebras in the theory of Hamiltonian equations. \emph{Jpn. J. Math.} 4 (2009), 141-252.

\bibitem{BFFLS}
F. Bayen, M. Flato, C. Fronsdal, A. Lichnerowicz and D. Sternheimer, Quantum mechanics as a deformation of classical mechanics. \emph{Lett. Math. Phys. }1 (1975/77), 521-530. 

\bibitem{BFFLS1}
F. Bayen, M. Flato, C. Fronsdal, A. Lichnerowicz and D. Sternheimer,  Deformation theory and quantization. I. Physical applications. \emph{Ann. Physics} 111 (1978), 61-110.


\bibitem{BFFLS2}
F. Bayen, M. Flato, C. Fronsdal, A. Lichnerowicz and D. Sternheimer,  Deformation theory and quantization. II. Physical applications. \emph{Ann. Physics} 111 (1978), 111-151.

\bibitem{CP1}V. Chari and A. Pressley, A Guide to Quantum Qroups. {Cambridge University Press}, Cambridge, 1994.

\bibitem{CK}
S. Cheng and V. G. Kac, Conformal modules. \emph{Asian J. Math.} 1 (1997), 181-193.

\bibitem{DAK}
A. D'Andrea and V. G. Kac, Structure theory of finite conformal algebras. \emph{Selecta Math. (N. S.)} 4 (1998), 377-418.

\bibitem{DK09}
A. De Sole and V. G. Kac, Lie conformal algebra cohomology and the variational complex. \emph{Comm. Math. Phys.} 292 (2009), 667-719.


\bibitem{Dolg}
I. A. Dolguntseva, The Hochschild cohomology for associative conformal algebras. \emph{Algebra Logic} 46 (2007), 373-384.

\bibitem{Dr87}
V. G. Drinfeld, Quantum groups. \emph{Proc. Internat. Congr. Math.}
(Berkeley, 1986), Amer. Math. Soc., Providence, RI, 1987, 798-820.


\bibitem{FGV} M. Flato, M. Gerstenhaber and A. A. Voronov, Cohomology and deformation of Leibniz pairs. \emph{Lett. Math. Phys.} 34 (1995), 77-90.

\bibitem{FB}
E. Frenkel and D. Ben-Zvi, Vertex Algebras and Algebraic Curves. Mathematical 
Surveys and Monographs, Vol. 88, Amer. Math. Soc., 2001.


\bibitem{Gel}
I. M. Gel'fand and I. Ya. Dorfman, Hamiltonian operators and algebraic structures associated with them.  \emph{Anal. i Prilozhen.} 13 (1979), 13-30.

\bibitem{GK04}
V. Ginzburg and D. Kaledin, Poisson deformations of symplectic quotient singularities.
\emph{Adv. Math.} 186 (2004), 1-57.


\bibitem{Kac97}
V.  G. Kac, Formal distribution algebras and conformal algebras, in: Brisbane Congress in Math. Phys., July 1997.

\bibitem{Kac}
V. G.  Kac, Vertex Algebras for Beginners. Second edition. University Lecture Series, 10. \emph{American Mathematical Society, Providence, RI}, 1998. vi+201 pp.


\bibitem{Kol06}
P. Kolesnikov, Associative conformal algebras with finite faithful representation. \emph{Adv.
	Math.} 202 (2006), 602-637.

\bibitem{Kol20}
P. S. Kolesnikov, Universal enveloping Poisson conformal algebras. \emph{Internat. J. Algebra Comput.} 30 (2020), 1015-1034.

\bibitem{Kol19}
P. S. Kolesnikov and R. A. Kozlov, On the Hochschild cohomologies of associative conformal algebras with a finite faithful representation. \emph{Commun. Math. Phys.} 369 (2019), 351-370.

\bibitem{KKP}
P. S. Kolesnikov, R. A. Kozlov and A. S. Panasenko,  Quadratic Lie conformal superalgebras related to Novikov superalgebras. \emph{J. Noncommut. Geom.} 15 (2021), 1485-1500. 

\bibitem{Kon03}
M. Kontsevich, Deformation quantization of Poisson manifolds.
\emph{Lett. Math. Phys.} 66 (2003), 157-216.

 \bibitem{Kubo1} F. Kubo, Finite-dimensional non-commutative Poisson algebras. \emph{J. Pure Appl. Algebra} 113 (1996), 307-314.

%\bibitem{Kubo2} F. Kubo, Noncommutative Poisson algebra structures on poset algebras and morphisms
%of Leibniz pairs. \emph{Bull. Kyushu Inst. Tech. Math. Natur. Sci.} 44 (1997), 1-5.

\bibitem{Kubo3} F. Kubo, Finite-dimensional simple Leibniz pairs and simple Poisson modules. \emph{Lett. Math. Phys.} 43 (1998), 21-29.

\bibitem{Kubo4} F. Kubo, Non-commutative Poisson algebra structures on affine Kac-Moody aglebras. \emph{J. Pure Appl. Algebra} 126 (1998), 267-286.


\bibitem{LHS}
H. Li, Vertex algebras and vertex Poisson algebras. \emph{Commun. Contemp. Math.} 6 (2004), 61-110.


\bibitem{Lic77}
A. Lichnerowicz, Les vari\'{e}t\'{e}s de Poisson et leur alg\`{e}bres de Lie associ\'{e}es. \emph{J. Diff. Geom.} 12 (1977), 253-300.


\bibitem{LZY}
J. Liu, S. Zhou and L. Yuan, Conformal $r$-matrix-Nijenhuis structures, symplectic-Nijenhuis structures and $\huaO N$-structures, accepted by J. Math. Phys., arXiv: 2204.11389. 


\bibitem{Pol97}
A. Polishchuk, Algebraic geometry of Poisson brackets. \emph{J. Math. Sci.}  84 (1997), 1413-1444.

\bibitem{Prim}
M. Primc, Vertex algebras generated by Lie algebras. \emph{J. Pure Appl. Algebra} 135 (1999), 253-293.

\bibitem{Rot1}
M. Roitman, On free conformal and vertex algebras. \emph{J. Algebra} 217 (1999), 496-527.

\bibitem{Rot2}
M. Roitman, Universal enveloping conformal algebras. \emph{Selecta Math. (N.S.)} 6 (2000), 319-345.



\bibitem{Vaisman1}I. Vaisman, Lectures on the Geometry of Poisson Manifolds. \emph{Progr. Math.} 118, Birkh\"{a}user Verlag, Basel, 1994.


\bibitem{Wei77}A. Weinstein, Lecture on Symplectic Manifolds. \emph{CBMS Regional Conference Series in Mathematics} 29, Amer. Math. Soc., Providence, R.I., 1979.

\bibitem{Xu}
P. Xu, Noncommutative Poisson algebras. {\em Amer. J. Math.} 116 (1994), 101-125.


\bibitem{XP}
X. Xu, Quadratic conformal superalgebras. \emph{J. Algebra} 231 (2000), 1-38.

\bibitem{YYY}
Y. Yang, Y. Yao and Y. Ye, (Quasi-)Poisson enveloping algebras. \emph{Acta Math. Sin.  (Engl. Ser.)} 29 (2013), 105-118.

\bibitem{Yuan}
L. Yuan, $\huaO$-operators and Nijenhius operators of associative conformal algebras. \emph{J. Algebra} 609 (2022), 245-291.






\end{thebibliography}
\end{document}